\documentclass[10pt]{article}
\usepackage{amssymb,amsmath}
\newcommand{\esp}{\quad}

\newcommand{\refe}[1]{(\ref{#1})}
\newcommand{\dst}{\displaystyle}
\newcommand{\NN}{{\mathbb N}}
\newcommand{\CC}{{\mathbb C}}
\newcommand{\tilG}{\widetilde{\Gamma}_q}
\newcommand{\fhyp}{\mbox{F}}
\renewcommand{\u}{\textrm{u}}
\newcommand{\bq}{\begin{equation}}
\newcommand{\eq}{\end{equation}}
\newcommand{\bt}{\begin{tabular}}
\newcommand{\et}{\end{tabular}}
\newcommand{\ba}{\begin{array}}
\newcommand{\ea}{\end{array}}
\newcommand{\qdiv}{\varkappa_q}

\newcommand{\btd}{\nabla}
\newcommand{\btu}{\Delta}

\newcommand{\half}{\mbox{$\frac{1}{2}$} }

\newcommand{\serie}{ \mbox{{\large$\varphi$}}}

\newcommand{\qnf}[1]{\mbox{$[{#1}]_q!$}}

\title{A $q$-analog of the Racah polynomials and
the $q$-algebra $SU_q(2)$\footnote{\noindent\textbf{Classification MSC numbers:} 33D80, 33D45
\hfill \tt arXiv:math.QA/0412540}
}
\date{November 25, 2004}
\author{R. \'Alvarez-Nodarse${}^\dag$, Yu. F. Smirnov${}^\S$, and
R. S. Costas-Santos${}^\star$ \\[5mm]
\small ${}^\dag$Departamento de An\'alisis Matem\'atico.\\
\small Universidad de Sevilla. Apdo. 1160, E-41080 Sevilla, Spain\\
\small${}^\S$Skobeltsyn Institute of Nuclear Physics.
Moscow State University.\\ \small Vorob'evy Gory, Moscow 119992, Russia \\
\small${}^\star$Departamento de Matem\'aticas.
E.P.S., Universidad Carlos III de Madrid.\\
\small Ave. Universidad 30, E-28911, Legan\'es, Spain
}

\begin{document}
\maketitle
\begin{abstract}
We study some $q$-analogues of the Racah polynomials and some
of their applications in the theory of representation of quantum
algebras.
\end{abstract}

\section{Introduction}

In the paper \cite{a-w-1} an orthogonal polynomial family that
generalizes the Racah coefficients or $6j$-symbols was introduced: the
so-called Racah and $q$-Racah polynomials. These polynomials were in the top of
the so-called Askey Scheme (see e.g. \cite{ks}) that contains all
classical families of hypergeometric orthogonal polynomials. Some
years later the same authors \cite{a-w-2} introduced the
celebrated Askey-Wilson polynomials. One of the important properties
of these polynomials is that from them one can obtain all known
families of hypergeometric polynomials and $q$-polynomials as
particular cases or as limit cases (for a review on this see the
nice survey \cite{ks}). The main tool in these two works was the
hypergeometric and basic series, respectively. On the other hand,
the authors of \cite{nu83} (see also \cite[Russian Edition]{nsu})
considered the $q$-polynomials as the solution of a second order
difference equation of hypergeometric-type on the non-linear
lattice $x(s)= c_1q^{s}+c_2q^{-s}+c_3$. In particular, they show
that the solution of the hypergeometric-type equation can be
expressed as certain basic series and, in such a way, they
recovered the results by Askey \& Wilson.

The interest of such polynomials increase after the appearance of the
$q$-algebras and quantum groups \cite{dri87,fad86,jim86,ku81,sk82}.
However, from the first attends to built the $q$-analog of the Wigner-Racah
formalism for the simplest quantum algebra $U_q(su(2))$ \cite{ki88}
(see also \cite{renato,ran-smi,mal92}) becomes clear that for obtaining
the $q$-polynomials intimately connected with the $q$-analogues of the Racah
and Clebsh Gordan coefficients, i.e., a $q$-analogue of the
Racah polynomials $u^{\alpha, \beta}_{n}(x(s),a,b)_q$ and the dual Hahn
polynomials $w^{c}_{n}(x(s),a,b)_q$, respectively, it is better to use a
different lattice ---in fact the $q$-Racah polynomials $R_n^{\beta,\gamma}(x(s),N,\delta)_q$
introduced in \cite{a-w-2} (see also \cite{ks}) were defined on the lattice $x(s)=q^{-s} +\delta q^{-N} q^{s}$
that depends not only of the variable $s$ but also on the parameters of the polynomials---,
namely,
\bq
x(s)=[s]_q [s + 1 ]_q,
\label{red-cua}
\eq
that only depends on $s$, where by $[s]_q$ we denote the $q$-numbers (in its symmetric form)
\bq\label{q-num}
[ s ]_q = \frac{q^{s/2} - q^{-s/2}}{q^{1/2} - q^{-1/2}},\quad
\forall s\in\CC.
\eq
With this choice the $q$-Racah polynomials
$u^{\alpha, \beta}_{n}(x(s),a,b)_q$ are proportional to the $q$-Racah
coefficients (or $6j$-symbols)  of the quantum algebra $U_q(su(2))$.
A very nice and simple approach to $6j$-symbols has been recently 
developed in \cite{ro03}. 

Moreover, this connection gives the possibility to a deeper study of the
Wigner-Racah formalism (or the $q$-analogue of the quantum theory of
angular momentum \cite{smi1,smi2,smi3,smi4}) for the quantum algebras
$U_q(su(2))$ and $U_q(su(1,1))$ using the powerful and well-known theory
of orthogonal polynomials on non-uniform lattices. On the other hand,
using the $q$-analogue of the quantum theory of
angular momentum \cite{smi1,smi2,smi3,smi4} we can obtain several results
for the $q$-polynomials, some of which are non trivial from the
point of view of the theory of orthogonal polynomials (see e.g.
the nice surveys \cite{koelink,vk}). In fact, in the present paper
we present a detailed study of some $q$-analogues of the Racah polynomials
on the lattice \refe{red-cua}: the $u^{\alpha, \beta}_{n}(x(s),a,b)_q$
and the $\widetilde{u}^{\alpha, \beta}_{n}(x(s),a,b)_q$ as well as
their connection with the $q$-Racah coefficients (or $6j$-symbols)
of the quantum algebra $U_q(su(2))$ in order to establish which properties
of the polynomials correspond to the $6j$-symbols and vice versa.

The structure of the paper is as follows: In section 2 we present
some general results from the theory of orthogonal polynomials on the
non-uniform lattices taken from \cite{ran,nsu}. In Section 2.1
a detailed discussion of the Racah polynomials
$u^{\alpha,\beta}_{n}(x(s),a,b)_q$ is presented, whereas in Section 2.2
the $\widetilde{u}^{\alpha, \beta}_{n}(x(s),a,b)_q$
are considered. In particular, a relation between these
families is established. In section 3 the comparative analysis of
such families and the $6j$-symbols of the quantum algebra $U_q(su(2))$
is developed which gives, on one hand, some information about the
Racah coefficients and, on the other hand, allow us to give a
group-theoretical interpretation of the Racah polynomials on the lattice
\refe{red-cua}. Finally, some comments and remarks about $q$-Racah
polynomials and the quantum algebra $U_q(su(3))$ are included.


\section{Some general properties of $q$-polynomials}

We will start with some general properties of orthogonal
hypergeometric polynomials on the non-uniform lattices \cite{arsuslov,nsu}.

The hypergeometric polynomials are the polynomial solutions $P_n(x(s))_q$ of the second order
linear difference equation of hypergeometric-type on the non-uniform lattice $x(s)$ (SODE)
\bq
\ba{c} \displaystyle
\sigma(s)\frac{\btu}{\btu x(s-\half)}  \frac{\btd y(s)}{\btd x(s)}
+ \tau(s)  \frac{\btu y(s)}{\btu x(s)}   + \lambda y(s) =0,\quad x(s)=c_1[q^{s}+q^{-s-\mu} ]+c_3,
\quad q^\mu=\frac{c_1}{c_2},
\\[5mm]
\btd f(s)=f(s)-f(s-1),\qquad\btu f(s)=f(s+1)-f(s),
 \ea
\label{eqdif-q}
\eq
or, equivalently
\bq
A_s y (s+1) + B_s y(s) + C_s y(s-1) + \lambda y(s) = 0,
\label{eqdif-q-new}
\eq
where
$$\ba{l} \dst
A_s= \dst\frac{\sigma (s) + \tau (s) \Delta x(s-\half)}
{\Delta x(s) \Delta x(s-\half)},
\quad
\dst C_s = \frac{\sigma (s) }{\nabla x(s) \Delta x(s-\half)},\quad
\dst B_s = -(A_s + C_s).
\ea
$$
Notice that $x(s)=x(-s-\mu)$.

In the following we will use the following notations\footnote{In the exponential
lattice $x(s)=c_1q^{\pm s}+c_3$, so $\mu=\pm\infty$, therefore
instead of using $\sigma(-s-\mu)$ one should use the equivalent
function $\sigma (s) + \tau (s) \Delta x(s-\half)$.} $P_n(s)_q:=P_n(x(s))_q$ and
$\sigma(-s-\mu)=\sigma (s) + \tau (s) \Delta x(s-\half)$. With this notation
the Eq. (\ref{eqdif-q}) becomes
\bq
\sigma(-s-\mu)\frac{\btu P_n(s)_q}{\btu x(s)} -
\sigma(s)\frac{\btd P_n(s)_q}{\btd x(s)} + \lambda_n \btu
x(s-\half) P_n(s)_q=0\,. \label{eqdif-q-equiv}
\eq
The polynomial solutions $P_n(s)_q$ of (\ref{eqdif-q}) can be
obtained by the following Rodrigues-type formula \cite{nsu,nu}
\bq
P_n(s)_q=\frac{B_n}{\rho(s)} \btd^{(n)} \rho _n(s),\qquad
\btd^{(n)}:= 
\frac{\btd}{\btd x_{1}(s)}\frac{\btd}{\btd  x_{2}(s)}
\cdots \frac{\btd}{\btd  x_n(s)},
\label{rodeq-q}
\eq
where $x_m(s)=x(s+\mbox{\footnotesize $\frac{m}{2}$})$,
\bq
\rho_n(s)=\rho(s+n) \prod_{m=1}^{n}\sigma(s+m),
\label{rhok-q}
\eq
and $\rho(s)$ is a solution of the Pearson-type equation
$\displaystyle{\Delta}  \left[\sigma(s) \rho(s)
\right]= \tau(s) \rho(s){\Delta x(s-1/2)}$, or equivalently,
\bq
\frac{\rho(s+1)}{\rho(s)}=  \frac{ \sigma(s)+\tau(s)\btu
x(s-\half)} {\sigma(s+1)}=\frac{ \sigma(-s-\mu)}{\sigma(s+1)}.
\label{pearson-q}
\eq
Let us point out that the function $\rho_n$ satisfy the equation
$\Delta \left[\sigma(s) \rho_n(s)\right]= \tau_n(s) \rho_n(s) \Delta x_n(s-1/2)$,
where $\tau_n(s)$ is given by
\bq
\tau_n(s) = \frac{\sigma(s+n)+\tau(s+n)\btu x(s+n-\frac{1}{2})-\sigma(s)}
{\btu x_{n-1}(s)}=\frac{\sigma(-s-n-\mu) -
\sigma(s)}{\btu x_n(s-\half) }=\tau_n' x_n(s)+\tau_n(0).
\label{tilde-tau-sig}
\eq
being
$$
\tau_n' =-\frac{\lambda_{2n+1}}{[2n+1]_q},\qquad \tau_n(0)=\frac{\sigma(-s_n^\star-n-\mu)-\sigma(s_n^\star)}
{x_n(s_n^\star+\half)-x_n(s_n^\star-\half)},
$$
where $s_n^\star$ is the zero of the function $x_n(s)$, i.e.,
$x_n(s_n^\star)=0$.

From \refe{rodeq-q} follows an explicit formula
for the polynomials $P_n$ \cite[Eq.(3.2.30)]{nsu}
\bq
P_n(s)_q=B_n\dst\sum_{m=0}^n\frac{[n]_q!(-1)^{m+n}}{[m]_q![n-m]_q!}
\frac{   \btd x(s+m-\mbox{\footnotesize $\frac{n-1}{2}$})  }
{  \displaystyle \prod _{l=0}^{n}   \btd x(s+\mbox{\footnotesize
$\frac{m-l+1}{2}$})   }
\frac{\rho_n(s-n+m)}{ \rho(s)},
\label{form-exp-q}
\eq
where $[n]_q$ denotes the {\em symmetric q-numbers} \refe{q-num}
and the $q$-factorials are given by
$$
[0]_q!:=1,\quad [n]_q!:=[1]_q[2]_q\cdots [n]_q,\quad n\in\NN.
$$
It can be shown \cite{arsuslov,nsu,nu}
that the most general polynomial solution of the
{\it q-}hypergeometric equation \refe{eqdif-q} corresponds to
\begin{equation}
\sigma(s) = A\prod_{i=1}^4 [s-s_i]_q = C q^{-2s}\prod_{i=1}^4(q^s-q^{s_i}),
\qquad A\cdot C\neq 0
\label{sigma-gen}
\end{equation}
and  has the form \cite[Eq. (49a), page 240]{nu}
\begin{equation}
\begin{array}{l}\displaystyle
P_n(s)_q= \displaystyle   D_n  \,\, {}_{4}
\mbox{\large{$\phi$}}_3 \left(\begin{array}{c} q^{-n},
q^{ 2\mu+n-1+\mbox{\tiny $\sum_{i=1}^4$} s_i}, q^{s_1-s},q^{s_1+s+\mu} \\
q^{s_1+s_2+\mu},q^{s_1+s_3+\mu}, q^{s_1+s_4+\mu} \end{array}
;\, q\, ,\, q \right),
\end{array}
\label{rep-serie-q}
\end{equation}
where the normalizing factor $D_n$ is given by ($\qdiv:=q^{1/2}-q^{-1/2}$)
$$
D_n\!=\dst\! B_n\!\left(\frac{-A}{c_1q^{\mu}\qdiv^5}\right)^n
\!\!\!q^{-\frac{n}{2}(3s_1+s_2+s_3+s_4+\frac{3(n-1)}{2})}
(q^{s_1+s_2+\mu};q)_n
(q^{s_1+s_3+\mu};q)_n (q^{s_1+s_4+\mu};q)_n.
$$
The basic hypergeometric
series ${}_r  \mbox{\large{$\phi$}}_p$ are defined by \cite{ks}
$$ \mbox{\small $\displaystyle
\begin{array}{l} \displaystyle _{r} \mbox{\large{$\phi$}}_p
\left(\!\begin{array}{c} {a_1,\dots,a_{r}} \\
 {b_1,\dots,b_p} \end{array} \,;\, q\, ,\, z \!\right)=
 \displaystyle \sum _{k=0}^{\infty}\frac{ (a_1;q)_k \cdots  (a_{r};q)_k}
{(b_1;q)_k  \cdots  (b_p;q)_k}\frac{z^k}{(q;q)_k}
\left[ (-1)^k q^{\frac k2(k-1)}\right]^{p-r+1} ,
\end{array} $}
$$
where $(a;q)_k= \prod_{m=0}^{k-1}(1-aq^m)$, is the  $q$-analogue of the Pocchammer symbol.

In this paper we will deal with  orthogonal {\it q-}polynomials
and functions. It can be proven \cite{nsu}, by using the
difference equation of hypergeometric-type \refe{eqdif-q}, that if
the boundary conditions $\sigma(s) \rho(s) x^{k}(s-{1}/{2})
\big|_{s=a,b} = 0$, for all $k \ge  0$, holds, then the
polynomials $P_n(s)_q$ are orthogonal with respect to the weight
function $\rho$, i.e.,
\begin{equation}
\begin{array}{c} \displaystyle
\sum_{s = a }^{b-1} P_n(s)_q P_m(s)_q\rho(s) \Delta
x(s-1/2) = \delta_{nm} d_n^2, \quad \!\!s=a,a+1,\dots,b-1.
\end{array}
\label{norm}
\end{equation}
The squared norm in \refe{norm} is
\cite[Eq. (3.7.15)]{nsu}
\bq
d_n^2 =(-1)^n A_{n,n} B_n^2 \sum_{s=a}^{b-n-1}\rho_n(s)
\Delta  x_n(s-1/2),
\label{norma}
\eq
where \cite[page 66]{nsu}
\bq\hspace{.1cm}
A_{n,k} = \frac{[n]_q!}{[n-k]_q!} \prod_{m=0}^{k-1}
\left(-\frac{\lambda_{n+m}}{[n+m]_q}\right) .
\label{A_nm-q}
\eq

A simple consequence of the orthogonality is the three-term recurrence relation (TTRR)
\bq
\label{ttrr}
 x(s)P_n(s)_q=\alpha_n P_{n+1}(s)_q+\beta_n P_n(s)_q+
\gamma_n P_{n-1}(s)_q,
\end{equation}
where $\alpha_n$, $\beta_n$ and $\gamma_n$ are given by
\bq
\alpha_n= \frac{a_n}{a_{n+1}} ,\quad
\beta_n = \frac{b_n}{a_{n}}-\frac{b_{n+1}}{a_{n+1}},\quad
\gamma_n= \frac{a_{n-1}}{a_{n}}\frac{d_{n}^2}{d_{n-1}^2},
\label{coef-RRTT-q}
\eq
being $a_n$ and $b_n$ the first and second coefficients in the power expansion
of $P_n$, i.e., $P_n(s)_q=a_n x^n(s)+b_n x^{n-1}(s)+ \cdots\,$.
Substituting $s=a$ in \refe{ttrr} we find
\bq
\beta_n=\frac{x(a)P_n(a)_q-\alpha_n P_{n+1}(a)_q-
\gamma_n P_{n-1}(a)_q}{P_n(a)_q},
\label{beta_n}
\eq
which is an alternative way for finding the coefficient $\beta_n$. Also we
can use the expression \cite[page 148]{ran}
$$
\beta_n=\frac{[n]_q\tau_{n-1}(0)}{\tau'_{n-1}}-
\frac{[n+1]_q\tau_{n}(0)}{\tau'_{n}}+c_3([n]_q+1-[n+1]_q).
$$
To compute $\alpha_n$ (and $\beta_n$) we need the following formulas (see e.g.
\cite[page 147]{ran})
\bq
a_n= \frac{B_n A_{n,n}}{[n]_q!},\qquad
\frac{b_n}{a_n}=\frac{[n]_q\tau_{n-1}(0)}
{\tau'_{n-1}}+c_3([n]_q-n).
\label{a_n-q}
\eq

The explicit expression of $\lambda_n$ is \cite[Eq. (52) page 232]{nu}
\bq\begin{split} \lambda_n
=&\dst-\frac{A q^{\mu}}{c_1^2(q^{1/2}-q^{-1/2})^4}[n]_q
\left[s_1+s_2+s_3+s_4+ 2\mu+n-1\right]_q\, \\
=&\dst -\frac{C\,q^{-n+1/2}}{c_1^2(q^{1/2}-q^{-1/2})^2 }\left(1-q^n\right)
\left(1-q^{s_1+s_2+s_3+s_4+2\mu+n-1}\right),
\end{split}
\label{lambda-q-equiv}
\eq
which can be obtained equating the largest powers of $q^s$ in
\refe{eqdif-q-equiv}.

From the Rodrigues formula  \cite[\S 5.6]{nsu-r,ran} follows that
\bq\ba{c}
\displaystyle
\frac{\btu P_n(s-\half)_q}{\btu x(s-\half)} = \frac{-\lambda_n B_n}
{\widetilde{B}_{n-1}}
\widetilde{P}_{n-1}(s)_q,
\ea
\label{for-dif}
\eq
where $\widetilde{P}_{n-1}$  denotes the polynomial orthogonal with
respect to the weight function $\widetilde{\rho}(s)= \rho_1(s-\half)$. On
the other hand, rewriting \refe{eqdif-q} as
$$
\left(\sigma(s)\frac{\btd}{\btd x_1(s)}+\tau(s) I\right)
\frac{\btu }{\btu x(s)} P_n(s)_q=-\lambda_n P_n(s)_q,
$$
it can be substituted by the following two first-order
difference equations
\bq
\frac{\btu }{\btu x(s)} P_n(s)_q= Q(s),\qquad
\left(\sigma(s)\frac{\btd}{\btd x_1(s)}+\tau(s) I\right) Q(s)=-\lambda_n
P_n(s)_q.
\label{difeq-1.1}
\eq
Using the fact that $\frac{\btu }{\btu x(s)} P_n(s)_q$ is a
polynomial of degree $n-1$ on $x(s+1/2)$ (see \cite[\S 3.1]{nsu})
it follows that
$$
\frac{\btu }{\btu x(s)} P_n(s)_q=C_n {Q}_{n-1}(s+\half),
$$
where $C_n$ is a normalizing constant. Comparison with \refe{for-dif}
implies that $Q(s)$ is the polynomial $\widetilde{P}_{n-1}$
orthogonal with respect to the function $\rho_1(s-\half)$
and $C_n={-\lambda_n B_n}/{\widetilde{B}_{n-1}}$. Therefore,
the second expression in \refe{difeq-1.1} becomes
\bq
P_n(s)_q= \frac{B_n}{\widetilde{B}_{n-1}}
\left(\sigma(s)\frac{\btd}{\btd x_1(s)}+\tau(s) I\right)
\widetilde{P}_{n-1}(s+\half)_q.
\label{dif-for-2}
\eq

The $q$-polynomials satisfy the following differentiation-type formula \cite[\S 5.6.1]{nsu-r,ran}
\bq
\sigma(s)\frac{\nabla P_n(s)_q}{\nabla
x(s)}=\frac{\lambda_n}{[n]_q\tau_n'}\left[\tau_n(s)P_n(s)_q-
\frac{B_n}{B_{n+1}}P_{n+1}(s)_q\right].
\label{d3a}
\eq
Then, using the explicit expression for the coefficient $\alpha_n$, we find
\begin{equation}
\label{d3} \sigma(s)\frac{\nabla P_n(s)_q}{\nabla
x(s)}=\frac{\lambda_n}{[n]_q} \frac{\tau_n(s)}{\tau_n'}P_n(s)_q-
\frac{\alpha_n\lambda_{2n}}{[2n]_q}P_{n+1}(s)_q .
\end{equation}

{}From the above equation using the identity
$\btu  \frac{\btd  P_{n}(s)_q}{\btd x(s)}=
\frac{\btu  P_{n}(s)_q}{\btu x(s)} - \frac{\btd  P_{n}(s)_q}{\btd x(s)}$
as well as the SODE  \refe{eqdif-q-equiv} we find
\bq \label{d4}
\begin{array}{c}\dst
\sigma(-s-\mu) \frac{\btu P_{n}(s)_q}{\btu x(s)}
=  \frac{\lambda_n}{[n]_q \tau'_n}\left[\big(\tau_n(s)-[n]_q
\tau'_n\btu x(s-\half)\big)P_{n}(s)_q-\frac{B_{n}}{B_{n+1}}
P_{n+1}(s)_q \right].
\end{array}
\end{equation}

To conclude this section we will introduce the following notation by
Nikiforov and Uvarov \cite{nsu,nu}. First we define another
{\em q-analog of the Pocchammer symbols} \cite[Eq. (3.11.1)]{nsu}
\bq
(a|q)_k = \prod_{m=0}^{k-1} [a+m]_q = \frac{\tilG(a+k)}{\tilG(a)}=
(-1)^k(q^a;q)_k (q^{1/2}-q^{-1/2})^{-k}q^{-\frac k4(k-1)-\frac{ka}2},
\label{sim-poc-q}
\eq
where $\tilG(x)$ is the $q$-analog of the $\Gamma$ function
introduced in \cite[Eq. (3.2.24)]{nsu}, and related to the classical
$q$-Gamma function $\Gamma_q$ by formula
$$
\tilG(s)= q^{-\frac{(s-1)(s-2)}{4}}\Gamma_q(s)=q^{-\frac{(s-1)(s-2)}{4}}
(1-q)^{1-s}\frac{(q;q)_\infty}{(q^s;q)_\infty},\quad 0<q<1.
$$

Next we define the {\em q-hypergeometric function}
${}_r\fhyp_p(\cdot|q,z)$
\bq \quad\quad \mbox{$\dst{}_{r}\fhyp_p
\left(\!\!\ba{c} {a_1,\dots,a_{r}} \\ {b_1,\dots,b_p} \ea \,
\bigg|\, q\,,\,  z \right)\!=\!
 \displaystyle
\sum _{k=0}^{\infty}\frac{ (a_1|q)_k(a_2|q)_k \cdot \cdot \cdot
(a_{r}|q)_k} {(b_1|q)_k(b_2|q)_k \cdot \cdot \cdot
(b_p|q)_k}\frac{z^k}{(1|q)_k} \left[\qdiv^{-k} q^{\frac{1}4
k(k-1)}\right]^{p-r+1} $}\!\!\!, \label{q-hip-def} \eq where, as
before, $\qdiv=q^{1/2}-q^{-1/2}$, and $(a|q)_k$ are given by
(\ref{sim-poc-q}). Notice that
$$
\lim_{q\to1}{}_{r}\fhyp_p\left(\ba{c}{a_1,a_2,\dots,a_{r}}\\
{b_1,b_2,\dots,b_p} \ea
 \,\bigg|\, q\,,\,  z \, \qdiv^{p-r+1} \right)=
\displaystyle \sum _{k=0}^{\infty}\frac{ (a_1)_k \cdots  (a_{r})_k}
{(b_1)_k \cdots  (b_p)_k}\frac{z^k}{k!} =
\displaystyle{}_{r}
\fhyp _p\displaystyle\bigg(\ba{c}{a_1,a_2,\dots,a_{r}}\\{b_1,b_2,\dots,b_p} \ea
\bigg| z\bigg),
$$
and
\bq
\ba{c}
{}_{p+1}\fhyp_p \left(\ba{c} {a_1,a_2,\dots,a_{p+1}} \\
{b_1,b_2,\dots,b_p} \ea \,\bigg|\, q\,,\, t \right)
\Bigg|_{\displaystyle t=t_0}=
\dst _{p+1}\serie_p
\left(\ba{c} {q^{a_1},q^{a_2},\dots,q^{a_{p+1}}} \\
{q^{b_1},q^{b_2},\dots,q^{b_p}}\ea \,\bigg|\, q\,,\, z \right),
\ea
\label{ser-q-fhyp}
\eq
where $t_0= z\,q^{\half \left(\sum_{i=1}^{p+1}a_i -\sum_{i=1}^{p} b_i -1 \right)}$.

Using the above notation the polynomial solutions of
\refe{eqdif-q} is \cite[Eq. (49), page 232]{nu}
\bq
\ba{l}\dst
P_n(s)_q=\dst B_n\left(\frac{A}{c_1q^{-\frac{\mu}{2}}
\qdiv^2} \right)^n(s_1+s_2+\mu|q)_n
(s_1+s_3+\mu|q)_n \times \\[7mm]\qquad\dst (s_1+s_4+\mu|q)_n \,\,
{}_{4}\fhyp_3 \left(\ba{c} -n, 2\mu+n-1+\dst \sum_{i=1}^4 s_i,
s_1-s,s_1+s+\mu \\ s_1+s_2+\mu,s_1+s_3+\mu, s_1+s_4+\mu \ea
\,\bigg|\, q\,,\, 1 \right).
\ea
\label{rep-hip-q}
\eq


\subsection{The $q$-Racah polynomials}

Here we will consider the $q$-Racah polynomials
$u_n^{\alpha,\beta}(x(s),a,b)_q$ on the lattice
$ x(s)=[s]_{q}[s+1]_{q} $ introduced in
\cite{ran-tesis,mal92,nsu}. For this lattice one has
\bq
c_1={q^{\half}}{\qdiv^{-2}},\quad
\mu=1,\quad c_3= -({q^{\half}+q^{-\half}}){\qdiv^{-2}}.
\label{red-cua-q}
\eq
Let chose $\sigma$ in \refe{sigma-gen} as
$$\dst
\sigma(s)=-\frac{q^{-2s}}{\qdiv^{4}q^{\frac{\alpha+\beta}2}}%
(q^s-q^a)(q^s-q^{-b})(q^s-q^{\beta-a})(q^s-q^{b+\alpha})=
[s-a]_{q} [s+b]_{q}[s+a-\beta]_{q}[b+\alpha-s ]_{q},
$$
i.e., $s_1=a,s_2=-b,s_3=\beta-a,s_4=b+\alpha$, $C=-q^{-\half(\alpha+\beta)}\qdiv^{-4}$, $A=-1$, and let
$B_n={(-1)^n}/{[n]_{q}!}$. Here, as before, $\qdiv=q^{1/2}-q^{-1/2}$.
Now from  (\ref{lambda-q-equiv}) we find
$$
\lambda_n=q^{-\half(\alpha+\beta+2n+1)}\qdiv^{-2}
(1-q^n)(1-q^{\alpha+\beta+n+1})=
[n]_q[n+\alpha+\beta+1]_q.
$$

To obtain $\tau_n(s)$ we use \refe{tilde-tau-sig}. In this case
$x_n(s)=[s+n/2]_q[s+n/2+1]_q$, then, choosing $s_n^\star=-n/2$, we get
\bq \tau_n(s)=\tau_n' x_n(s)+\tau_n(0),\quad
\tau_n'=-[2n+\alpha + \beta+2]_q,\quad
\tau_n(0)=\sigma(-n/2-1)-\sigma(-n/2). \label{tau_n-Racah} \eq
Taking into account that  $\tau(s)=\tau_0(s)$, we obtain the
corresponding function $\tau(s)$
$$
\tau(s)=-[2+\alpha + \beta]_q x(s)+\sigma(-1)-\sigma(0).
$$

\subsubsection{The orthogonality and the norm $d_n^2$}

A solution of the Pearson-type difference equation \refe{pearson-q} is
$$
\rho(s)= \frac{\tilG(s+a+1) \tilG(s-a+\beta+1)\tilG(s+\alpha+b+1)
\tilG(b+\alpha-s)}{ \tilG(s-a+1) \tilG(s+b+1)\tilG(s+a-\beta+1)
\tilG(b-s)}.
$$
Since $\sigma(a)\rho(a)=\sigma(b)\rho(b)=0$, then
the $q$-Racah polynomials satisfy the orthogonality
relation
$$
\sum_{s=a}^{b-1} u_n^{\alpha,\beta}(x(s),a,b)_q
u_m^{\alpha,\beta}(x(s),a,b)_q \rho(s) [2s+1]_q = 0,\esp \esp n\neq m,
$$
with the restrictions $-\half<a\leq b-1, \alpha>-1, -1<\beta<2a+1$.
Let us now compute the square of the norm $d_n^2$. From \refe{rhok-q} and \refe{A_nm-q} follow
$$\begin{array}{l}
\rho_n(s)\!=\dst\frac{\widetilde{\Gamma}_q(s+n+a+1)
\widetilde{\Gamma}_q(s+n-a+\beta+1)
\widetilde{\Gamma}_q(s+n+\alpha+b+1)
\widetilde{\Gamma}_q(b+\alpha-s)}
{\widetilde{\Gamma}_q(s-a+1)\widetilde{\Gamma}_q(s+b+1)
\widetilde{\Gamma}_q(s+a-\beta+1)
\widetilde{\Gamma}_q(b-s-n)},\!\!
\end{array}
$$
$$
A_{n,n}=[n]_q! (-1)^n
\frac {\widetilde{\Gamma}_q(\alpha+\beta+2n+1)}
{\widetilde{\Gamma}_q(\alpha+\beta+n+1)}\quad\Rightarrow\quad
\Lambda_n:=
(-1)^n A_{n,n} B_n^2=\frac {\widetilde{\Gamma}_q(\alpha+\beta+2n+1)}
{[n]_q!\widetilde{\Gamma}_q(\alpha+\beta+n+1)}.
$$
Taking into account that $\dst \nabla x_{n+1}(s)=[2s+n+1]_q$, using \refe{norma},
and the identity
\bq
\tilG(A-s)=\frac{\tilG(A)(-1)^s}{(1-A|q)_s},
\label{ide-gam}
\eq
we have
$$
\begin{array}{l}
\dst
d_n^2=\!\Lambda_n\!\!
\sum_{s=a}^{b-n-1}\!\frac{\widetilde{\Gamma}_q(s+n+a+\!1)
\widetilde{\Gamma}_q(s+n-a+\beta+\!1)
\widetilde{\Gamma}_q(s+n+\alpha+b+\!1)
\widetilde{\Gamma}_q(b+\alpha-s)}
{\widetilde{\Gamma}_q(s-a+1)\widetilde{\Gamma}_q(s+b+1)
\widetilde{\Gamma}_q(s+a-\beta+1)
\widetilde{\Gamma}_q(b-s-n)[2s+n+1]_q^{-1}}
\\[0.8cm]
\dst =
\!\Lambda_n\!\!
\sum_{s=0}^{b-a-n-1}\!\frac{\widetilde{\Gamma}_q(s+n+2a+\!1)
\widetilde{\Gamma}_q(s+n+\beta+\!1)
\widetilde{\Gamma}_q(s+n+\alpha+b+a+\!1)
\widetilde{\Gamma}_q(b-a+\alpha-s)}
{\widetilde{\Gamma}_q(s+1)\widetilde{\Gamma}_q(s+b+a+1)
\widetilde{\Gamma}_q(s+2a-\beta+1)
\widetilde{\Gamma}_q(b-a-s-n)[2s+2a+n+1]_q^{-1}}\\[8mm]
\dst=
\frac{\widetilde{\Gamma}_q(\alpha+\beta+2n+1)
\widetilde{\Gamma}_q(2a+n+1)
\widetilde{\Gamma}_q(n+\beta+1)\widetilde{\Gamma}_q(a+b+n+\alpha+1)
\widetilde{\Gamma}_q(b+\alpha-a)}{[n]_q!\widetilde{\Gamma}_q(\alpha+\beta+n+1)
\widetilde{\Gamma}_q(a+b+1)\widetilde{\Gamma}_q(2a-\beta+1)\widetilde{\Gamma}_q(b-a-n)}\times
\\[0.8cm]\dst\quad  \sum_{s=0}^{b-a-n-1}
\frac{(n+2a+1,n+\beta+1,n+a+\alpha+b+1,1-b+a+n|q)_{s}}
{(1,a+b+1,2a-\beta+1,1-b+a-\alpha|q)_{s}}
[2s+2a+n+1]_q.
\ea
$$
In the following we denote by $S_n$ the sum in the last
expression. If we now use that
$\dst(a|q)_n={(-1)^n(q^a;q)_n}{q^{-\frac{n}{4}(n+2a-1)}\qdiv^{-n}}$, as well as
the identity
$$
[2s+2a+n+1]_q=q^{-s}[2a+n+1]_q
\frac{(q^{a+\frac{n+1}2+1};q)(-q^{a+\frac{n+1}2+1};q)_s}
{(q^{a+\frac{n+1}2};q)(-q^{a+\frac{n+1}2};q)},
$$
we obtain
\begin{small}
$$
\ba{l}
S_n=\!\!\!\dst \sum_{s=0}^{b-a-n-1}
\!\frac{(q^{2a+n+\!1},q^{n+\beta+1},q^{n+\alpha+b+a+\!1},
q^{1-b+a+n},q^{\half (2a+n+\!3)},
-q^{\half (2a+n+\!3)};q)_s}{(q,q^{a+b+1},q^{2a-\beta+1},
q^{1-b-\alpha+a},q^{\half (2a+n+1)},-q^{\half (2a+n+1)};q)_s[2a\!+\!n\!+\!1]_q^{-1}}
q^{-s(1+2n+\beta+\alpha)}
\\[0.8cm]
=[2a\!+\!n\!+\!1]_q \,_6 \phi_5 \!\! \left( \!\!\!
\begin{array}{c}
q^{2a+n+1}\!,q^{n+\beta+1},q^{n+a+\alpha+b+1}\!,q^{1-b+a+n}\!,q^{\half (2a+n+\!3)}\!,
-q^{\half (2a+n+\!3)} \\
q^{a+b+1},q^{2a-\beta+1},q^{1-b-\alpha+a},q^{\half (2a+n+1)},-q^{\half (2a+n+1)}
\end{array}
\!\! \Bigg| \! \ q,q^{-1-2n-\beta-\alpha}
\right).
\ea
$$
\end{small}
But the above ${}_6\phi_5$ series is a very-well-poised ${}_6\phi_5$ basic series
and therefore by using the   summation formula
\cite[Eq. (II.21) page 238]{gasper1}
$$
_{6}\varphi_{5}\Bigg({a, \ qa^{1/2}, \ -qa^{1/2}, \ b, \ c, \ q^{-k}
\atop \ a^{1/2}, \ -a^{1/2}, aq/b, \ aq/c, aq^{k+1}}\Bigg| \ q,{aq^{k+1}\over bc}\Bigg) =
{(aq, aq/bc;q)_{k}\over (aq/b,aq/c;q)_{k}},
$$
with $k=b-a-n-1$, $a=q^{2a+n+1}$, $b=q^{n+\beta+1}$,
$c=q^{n+a+\alpha+b+1}$, we obtain
$$\ba{l}\dst
S_n= [2a\!+\!n\!+\!1]_q
\frac{(q^{2a+n+2},q^{-n+a-b-\alpha-\beta};q)_{b-a-n-1}}
{(q^{2a-\beta+1},q^{a-b-\alpha+1};q)_{b-a-n-1}}\\[5mm]
\quad=\dst
[2a\!+\!n\!+\!1]_q
\frac{({2a+n+2}|q)_{b-a-n-1}({-n+a-b-\alpha-\beta}|q)_{b-a-n-1}}
{({2a-\beta+1}|q)_{b-a-n-1}({a-b-\alpha+1}|q)_{b-a-n-1}}.
\end{array}
$$
Finally, using \refe{ide-gam} and \refe{sim-poc-q}
$$
S_n= [2a\!+\!n\!+\!1]_q\frac{\tilG(a+b+1)\tilG(2a-\beta+1)
\tilG(b-a+\alpha+\beta+n+1)\tilG(\alpha+n+1)}
{\tilG(n+2a+2)\tilG(b+a-\beta-n)\tilG(\alpha+\beta+2n+2)
\tilG(b-a+\alpha)},
$$
thus
$$\ba{rl}
\dst d_n^2=&\dst\frac{\widetilde{\Gamma}_q(\alpha+\beta+2n+1)
\widetilde{\Gamma}_q(2a+n+1)
\widetilde{\Gamma}_q(n+\beta+1)\widetilde{\Gamma}_q(a+b+n+\alpha+1)
\widetilde{\Gamma}_q(b+\alpha-a)}{[n]_q!\widetilde{\Gamma}_q(\alpha+\beta+n+1)
\widetilde{\Gamma}_q(a+b+1)\widetilde{\Gamma}_q(2a-\beta+1)
\widetilde{\Gamma}_q(b-a-n)} S_n\\[5mm]=& \dst
\frac{\tilG(\alpha+n+1)\tilG( \beta+n+1) \tilG(b-a+\alpha+\beta+n+1)
 \tilG(a+b+\alpha+n+1)}{[\alpha+\beta+2n+1]_q \tilG(n+1)
\tilG(\alpha+\beta+n+1) \tilG(b-a-n) \tilG(a+b-\beta-n)}.
\ea
$$

\subsubsection{The hypergeometric representation}

From formula \refe{rep-serie-q} and \refe{rep-hip-q}
the following two equivalent hypergeometric representations hold
\bq
\begin{split}\dst
u_n^{\alpha,\beta}(x(s),a,b)_q= & \dst
\frac{q^{-\frac{n}{2}(2a+\alpha+\beta+n+1)}
(q^{a-b+1};q)_n (q^{\beta+1};q)_n
 (q^{a+b+\alpha+1};q)_n}{\qdiv^{2n} (q;q)_n}\times\\  &
{}_{4}\serie_3 \left(\ba{c} q^{-n},q^{ \alpha+\beta+n+1}, q^{a-s},
q^{a+s+1}  \\ q^{a-b+1},q^{\beta+1},q^{ a+b+\alpha+1}  \ea
\,\bigg|\, q \,,\, q \right),
\end{split}
\label{pol-rac-nu}
\eq
and
\bq\label{pol-rac-nu-F}
\begin{split}\dst
u_n^{\alpha,\beta}(x(s),a,b)_q= & \dst \frac{(a-b+1|q)_n (\beta+1|q)_n
(a+b+\alpha+1|q)_n}{[n]_{q}!}\times\\  &
_{4}\fhyp_3 \left(\ba{c} -n, \alpha+\beta+n+1, a-s,a+s+1 \\
 a-b+1,\beta+1, a+b+\alpha+1 \ea \,\bigg|\, q \,,\, 1 \right).
\end{split}
\eq
Using the Sears transformation formula \cite[Eq. (III.15)]{gasper1}
we obtain the equivalent formulas
\bq
\begin{split}\dst
u_n^{\alpha,\beta}(x(s),a,b)_q= & \dst
\frac{q^{-\frac{n}{2}(-2b+\alpha+\beta+n+1)}
(q^{a-b+1};q)_n (q^{\alpha+1};q)_n(q^{\beta-a-b+1};q)_n}{\qdiv^{2n}
(q;q)_n}\times\\ & {}_{4}\serie_3
\left(\ba{c} q^{-n},q^{ \alpha+\beta+n+1}, q^{-b-s},
q^{-b+s+1}  \\ q^{a-b+1},q^{\alpha+1},q^{-a-b+\beta+1}\ea\,\bigg|\, q
\,,\, q \right),
\end{split}
\label{pol-rac-nu2}
\eq
and
\bq\label{pol-rac-nu-F-2}
\begin{split}\dst
u_n^{\alpha,\beta}(x(s),a,b)_q= & \dst \frac{(a-b+1|q)_n (\alpha+1|q)_n
(-a-b+\beta+1|q)_n}{[n]_{q}!}\times\\  &
_{4}\fhyp_3 \left(\ba{c} -n, \alpha+\beta+n+1, -b-s,-b+s+1 \\
 a-b+1,\alpha+1, -a-b+\beta+1 \ea \,\bigg|\, q \,,\, 1 \right).
\end{split}
\eq
\noindent\textbf{Remark:} From the above formulas follow that the polynomials $u_n^{\alpha,\beta}(x(s),a,b)_q$ are 
multiples of the standard $q$-Racah polynomials $R_n(\mu(q^{b+s});q^\alpha,q^\beta,q^{a-b},q^{-a-b}|q)$.\\

From the above hypergeometric representations also follow the values
\bq\label{val-a}\begin{split}
u_n^{\alpha,\beta}(x(a),a,b)_q=&\dst
 \frac{(a-b+1|q)_n (\beta+1|q)_n (a+b+\alpha+1|q)_n}{[n]_{q}!}\\
=&\dst \frac{(q^{a-b+1};q)_n (q^{\beta+1};q)_n
 (q^{a+b+\alpha+1};q)_n}{q^{\frac{n}{2}(2a+\alpha+\beta+n+1)}\qdiv^{2n} (q;q)_n},
\end{split}
\eq
\bq\label{val-b}
\begin{split}
u_n^{\alpha,\beta}(x(b-1),a,b)_q =&\dst  \frac{(a-b+1|q)_n
(\alpha+1|q)_n (-a-b+\beta+1|q)_n}{[n]_{q}!}\\
=&\dst \frac{(q^{a-b+1};q)_n (q^{\alpha+1};q)_n(q^{\beta-a-b+1};q)_n}
{q^{\frac{n}{2}(-2b+\alpha+\beta+n+1)}\qdiv^{2n} (q;q)_n}.
\end{split}
\eq
The formula \refe{form-exp-q} leads to the
following explicit formula\footnote{Obviously the formulas
\refe{pol-rac-nu} and \refe{pol-rac-nu2} also give equivalent
explicit formulas.}
\bq\label{exp-for-rac}
\begin{split}
u_n^{\alpha,\beta}(x(s),a,b)_q  =&  \dst
\frac{ \tilG(s-a+1) \tilG(s+b+1)\tilG(s+a-\beta+1) \tilG(b-s) }
{\tilG(s+a+1) \tilG(s-a+\beta+1)\tilG(s+\alpha+b+1) \tilG(b+\alpha-s)}\times\\
&\dst \sum_{k=0}^{n} \frac{(-1)^k [2s+2k-n+1]_q
\tilG(s+k+a+1) \tilG(2s+k-n+1)}
{\tilG(k+1) \tilG(n-k+1) \tilG(2s+k+2)\tilG(s-n+k-a+1)}\times\\
&\quad\dst\frac{\tilG(s+k-a+\beta+1)\tilG(s+k+\alpha+b+1)
\tilG(b+\alpha-s+n-k)}
{\tilG(s-n+k+b+1)\tilG(s-n+k+a-\beta+1) \tilG(b-s-k)},
\end{split}
\eq
from where follows
\bq\ba{l}\dst
u_n^{\alpha,\beta}(x(a),a,b)_q =
\frac{(-1)^n \tilG(b-a)\tilG(\beta+n+1)\tilG(b+a+\alpha+n+1)}
{[n]!\tilG(b-a-n)\tilG(\beta+1)\tilG(b+a+\alpha+1)},\\[4mm]
\dst
u_n^{\alpha,\beta}(x(b-1),a,b)_q =
\frac{\tilG(b-a)\tilG(\alpha+n+1)\tilG(b+a-\beta)}
{[n]!\tilG(b-a-n)\tilG(\alpha+1)\tilG(b+a-\beta-n)},
\ea
\label{val-ext}
\eq
that coincide with the values \refe{val-a} and \refe{val-b}
obtained before.

{From} the hypergeometric representation the following symmetry property follows
$$
u_n^{\alpha,\beta}(x(s),a,b)_q=u_n^{-b-a+\beta,b+a+\alpha}
(x(s),a,b)_q.
$$

Finally, notice that from \refe{pol-rac-nu} (or \refe{pol-rac-nu2})
follows that $u_{n}^{\alpha,\beta}(x(s),a,b)_q$ is a polynomial of degree $n$
on $x(s)=[s]_q[s+1]_q$. In fact,
$$
(q^{a-s};q)_k (q^{a+s+1};q)_k=(-1)^k q^{k(a+\frac{k+1}2)}
\prod_{l=0}^{k-1}\left(\frac{x(s)-c_3}{c_1}-q^{-\frac12}
(q^{a+l+\frac12}+q^{-a-l-\frac12})\right),
$$
where $c_1$ and $c_3$ are given in \refe{red-cua-q}.

\subsubsection{Three-term recurrence relation and differentiation formulas}

To derive the coefficients of the TTRR \refe{ttrr} we use \refe{coef-RRTT-q} and
\refe{beta_n}. Using \refe{a_n-q} and \refe{coef-RRTT-q}, we obtain
$$
a_n=\frac{\widetilde{\Gamma}_q(\alpha+\beta+2n+1)}
{[n]_q!\widetilde{\Gamma}_q(\alpha+\beta+n+1)},\quad
\alpha_n=\frac{ [n+1]_q [\alpha+\beta+n+1]_q}
{[\alpha+\beta+2n+1]_q[\alpha+\beta+2n+2]_q}.
$$
To find $\gamma_n$ we use  \refe{coef-RRTT-q}
$$
\gamma_n=
\frac{[a+b+\alpha+n]_q[a+b-\beta-n]_q[\alpha+n]_q[\beta+n]_q
[b-a+\alpha+\beta+n]_q[b-a-n]_q}
{[\alpha+\beta+2n]_q[\alpha+\beta+2n+1]_q}.
$$
To compute $\beta_n$ we use \refe{beta_n}
$$\ba{rl}
\beta_n=&\dst x(a)-\alpha_n\frac{u_{n+1}^{\alpha,\beta}(x(a),a,b)_q}
{u_n^{\alpha,\beta}(x(a),a,b)_q}-
\gamma_n\frac{u_{n-1}^{\alpha,\beta}(x(a),a,b)_q}
{u_n^{\alpha,\beta}(x(a),a,b)_q}\\[5mm]=&\dst  [a]_q[a+1]_q -
\frac{[\alpha+\beta+n+1]_q[a-b+n+1]_q[\beta+n+1]_q[a+b+\alpha+n+1]_q}
{[\alpha+\beta+2n+1]_q[\alpha+\beta+2n+2]_q} \\[0.5cm]
\dst&\dst +\frac{[\alpha+n]_q[b-a+\alpha+\beta+n]_q[a+b-\beta-n]_q[n]_q}
{[\alpha+\beta+2n]_q[\alpha+\beta+2n+1]_q}.
\ea
$$

\begin{table}[t]
\caption{Main data of the {\em q-}Racah
polynomials $u_n^{\alpha,\beta}(x(s),a,b)_q $\label{tabla1}}

\vspace*{.5cm}

{\small
\begin{center}
{ \renewcommand{\arraystretch}{.35}
\begin{tabular}{|@{}c@{}| | @{}c@{}|}\hline
 &  \\
$P_n(s)$ & $u_n^{\alpha,\beta}(x(s),a,b)_q   \,,\esp x(s) = [s]_q [s+1]_q$ \\
 &  \\
\hline\hline
 &  \\
$ (a,b)$ & \mbox{
$[a,b-1]$} \\
 &  \\
\hline
 &  \\
$\rho(s)$ & $  \dst\frac{\tilG(s+a+1) \tilG(s-a+\beta+1)\tilG(s+\alpha+b+1)
 \tilG(b+\alpha-s)}{ \tilG(s-a+1) \tilG(s+b+1)\tilG(s+a-\beta+1) \tilG(b-s)}  $ \\
 &  \\
 &  \mbox{
$-\half<a\leq b-1, \alpha>-1,-1<\beta<2a+1$} \\
 &  \\
\hline
 &  \\
$\sigma(s)$ & \mbox{
$  [s-a]_{q} [s+b]_{q}[s+a-\beta]_{q}[b+\alpha-s ]_{q} $} \\
 &  \\
\hline
 &  \\
$\sigma(-s-1)$ & \mbox{
$  [s+a+1]_{q} [b-s-1]_{q}[s-a+\beta+1]_{q}[b+\alpha+s+1 ]_{q} $} \\
 &  \\
\hline
 &  \\
$\tau(s) $ &  \mbox{
$ [\alpha+1]_q [a]_q [a-\beta]_q  +
 [\beta+1]_q [b]_q   [b+\alpha]_q -[\alpha+1]_q  [\beta+1]_q -
[\alpha+\beta+2]_q x(s) $ } \\
 &  \\
\hline
 &  \\
$\tau_n(s) $ &  \mbox{
\begin{tabular}{r}
$ -[\alpha+\beta+2n+2]_q x(s+\frac n2)+[a+\frac n2+1]_{q}[b-\frac n2-1]_{q}
[\beta+\frac n2+1-a]_{q}[b+\alpha+\frac n2+1]_{q}$\\[3mm] $-[a+\frac n2]_{q}
[b-\frac n2]_{q}[\beta+\frac n2-a]_{q}[b+\alpha+\frac n2]_{q}$ \end{tabular}}\\
 &  \\
 \hline
 &  \\
$\lambda _n$ &\mbox{
$[n]_q[\alpha+\beta+n+1]_q$} \\
 &  \\
\hline
 &  \\
$B_n$ &  $  \dst\frac{(-1)^n}{[n]_q!}$
 \\
 &  \\
 \hline
 &  \\
$d_n^2$ & $  \dst\frac{\tilG(\alpha+n+1)\tilG(\beta+n+1)\tilG(b-a+\alpha+\beta+n+1)
 \tilG(a+b+\alpha+n+1)}{[\alpha+\beta+2n+1]_q\tilG(n+1)\tilG(\alpha+\beta+n+1)
 \tilG(b-a-n) \tilG(a+b-\beta-n)\, }\,\,$\, \\
 &  \\
 \hline
 &  \\
$\,\rho_n(s)\,$ & $ \dst\frac{\tilG(s+n+a+1) \tilG(s+n-a+\beta+1)
\tilG(s+n+\alpha+b+1) \tilG(b+\alpha-s)}
{ \tilG(s-a+1) \tilG(s+b+1)\tilG(s+a-\beta+1) \tilG(b-s-n)}  $  \\
 &  \\
 \hline
  & \\
$a_n$ & $ \dst\frac{\tilG[\alpha+\beta+2n+1]_q}{ [n]_q!\tilG[\alpha+\beta+n+1]_q}$
\\
 & \\
\hline
 &  \\
$\alpha _n$ &  $  \dst\frac{ [n+1]_q [\alpha+\beta+n+1]_q}
{[\alpha+\beta+2n+1]_q[\alpha+\beta+2n+2]_q}$ \\
 &  \\
\hline
 &  \\
 $\beta_n $ &  \begin{tabular}{c}
$\dst[a]_q[a+1]_q -\frac{[\alpha+\beta+n+1]_q[a-b+n+1]_q[\beta+n+1]_q[a+b+\alpha+n+1]_q}
{[\alpha+\beta+2n+1]_q[\alpha+\beta+2n+2]_q}$ \\[3mm]
$\dst+\frac{[\alpha+n]_q[b-a+\alpha+\beta+n]_q[a+b-\beta-n]_q[n]_q}
{[\alpha+\beta+2n]_q[\alpha+\beta+2n+1]_q}$
\end{tabular}\\
 &  \\
\hline
 &  \\
$\gamma_n$ & $ \dst\frac{[a+b+\alpha+n]_q[a+b-\beta-n]_q[\alpha+n]_q[\beta+n]_q
[b-a+\alpha+\beta+n]_q[b-a-n]_q}{[\alpha+\beta+2n]_q[\alpha+\beta+2n+1]_q} $ \\
 &  \\
\hline
\end{tabular}  }
\end{center}  }
\end{table}

The differentiation formulas \refe{for-dif} and \refe{dif-for-2} yield
\bq\label{dif-for-rac1}
\frac{\btu u_n^{\alpha,\beta}(x(s),a,b)_q }{\btu x(s)}= [\alpha+\beta+n+1]_q
u_{n-1}^{\alpha+1,\beta+1}(x(s+\half),a+\half,b-\half)_q,
\eq
\bq\label{dif-for-rac2}\begin{split}
-[n]_q[2s+1]_q u_n^{\alpha,\beta}(x(s),a,b)_q = &
\sigma(-s-1)u_{n-1}^{\alpha+1,\beta+1}(x(s+\half),a+\half,b-\half)_q\\[4mm] & -
\sigma(s) u_{n-1}^{\alpha+1,\beta+1}(x(s-\half),a+\half,b-\half)_q,
\end{split}\eq
respectively. Finally, formulas \refe{d3a} (or \refe{d3})
and \refe{d4} lead to the {\it differentiation formulas}
\bq \label{first-diff-form}
\sigma(s)\frac{\nabla u_n^{\alpha,\beta}(x(s),a,b)_q }{[2s]_q}=
-\frac{[\alpha+\beta+n+1]_q}{[\alpha+\beta+2n+2]_q}\Bigg[
\tau_n(s)u_n^{\alpha,\beta}(x(s),a,b)_q +[n+1]_q
u_{n+1}^{\alpha,\beta}(x(s),a,b)_q \Bigg],
\eq
\bq \label{scnd-diff-form}
\ba{l}\dst
\sigma(-s-1)\frac{\Delta u_{n}^{\alpha,\beta}(x(s),a,b)_q}{[2s+2]_q}=
-\frac{[\alpha+\beta+n+1]_q}{[\alpha+\beta+2n+2]_q}\times\\[5mm]
\dst \quad \Bigg[(\tau_n(s)+[n]_q[\alpha+\beta+2n+2]_q
[2s+1]_q)u_{n}^{\alpha,\beta}(x(s),a,b)_q +[n+1]_q u_{n+1}^{\alpha,\beta}(x(s),a,b)_q
\Bigg],
\ea
\eq
where $\tau_n(s)$ is given in \refe{tau_n-Racah}.
\subsubsection{The duality of the Racah polynomials}
\label{sec-dual-u}
In this section we will discuss the duality property of the
$q$-Racah polynomials $u_{n}^{\alpha,\beta}(x(s),a,b)_q$.
We will follow \cite[pages 38-39]{nsu}. First of all, notice that the orthogonal
relation \refe{norm} for the Racah polynomials can be written in the form
$$
\sum_{t=0}^{N-1} C_{tn} C_{tm}=\delta_{n,m},\quad
C_{tn}= \frac{u_{n}^{\alpha,\beta}(x(t+a),a,b)_q\sqrt{\rho(t+a)\btu x(t+a-1/2)}}{d_n},\quad
N=b-a,
$$
where $\rho(s)$ and $d_n$ are the weight function and the norm of the
$q$-Racah polynomials $u_{n}^{\alpha,\beta}(x(s),a,b)_q$, respectively. The above relation
can be understood as the orthogonality property of the orthogonal matrix $C=(C_{tn})_{t,n=0}^{N-1}$
by its first index. If we now use the orthogonality of $C$ by the second
index we get
$$
\sum_{n=0}^{N-1} C_{tn} C_{t'n}=\delta_{t,t'}, \quad N=b-a,
$$
that leads to the dual orthogonality relation for the $q$-Racah polynomials
\bq\label{dual-ort-u}
\sum_{n=0}^{N-1} u_{n}^{\alpha,\beta}(x(s),a,b)_q u_{n}^{\alpha,\beta}(x(s'),a,b)_q\frac1{d_n^2}
=\frac1{\rho(s)\btu x(s-1/2)}\delta_{s,s'}.
\eq

The next step is to identify the functions $u_{n}^{\alpha,\beta}(x(s),a,b)_q$
as polynomials on some lattice $x(n)$. Before starting let us mention that from
the representation  \refe{pol-rac-nu} and the identity
$$
(q^{-n};q)_k (q^{\alpha+\beta+n+1};q)_k=
\prod_{l=0}^{k-1}\left(1+ q^{\alpha+\beta+2l+1}-
 q^{\frac{\alpha+\beta+1}2+l}
\left( \qdiv^2 x(t)+q^{\frac12}+q^{-\frac12}\right)\right),
$$
where $x(t)=[t]_q[t+1]_q=[n+\frac{\alpha+\beta}2]_q[n+\frac{\alpha+\beta}2+1]_q$,
follows that $u_{n}^{\alpha,\beta}(x(s),a,b)_q$ also constitutes a polynomial
of degree $s-a$ (for $s=a,a+1,\ldots, b-a-1$) on $x(t)$ with $t=n+\frac{\alpha+\beta}2$.

Let us now define the polynomials ---compare with the definition
of the Racah polynomials \refe{pol-rac-nu-F}---
\bq\begin{split}
\u_k^{\alpha',\beta'}(x(t),a',b')_q=&
\frac{(-1)^{k} \tilG(b'-a')\tilG(\beta'+k+1)\tilG(b'+a'+\alpha'+k+1)}
{[k]!\tilG(b'-a'-k)\tilG(\beta'+1)\tilG(b'+a'+\alpha'+1)}
\times\\[4mm] &
_{4}\fhyp_3 \left(\ba{c} -k, \alpha'+\beta'+k+1, a'-t,a'+t+1 \\
 a'-b'+1,\beta'+1, a'+b'+\alpha'+1 \ea \,\bigg|\, q \,,\, 1 \right),
\end{split}
\label{dua-rac}
\eq
where
\bq\label{change}
k=s-a, \quad t=n+\frac{\alpha+\beta}2,\quad
a'=\frac{\alpha+\beta}2, \quad b'=b-a+\frac{\alpha+\beta}2,
\quad \alpha'=2a-\beta, \quad \beta'=\beta.
\eq
Obviously they are polynomials of degree $k=s-a$ on the lattice $x(t)$ that satisfy the
orthogonality property
\bq\label{rel-dual-ort-u}
\sum_{t=a'}^{b'-1} \u_k^{\alpha',\beta'}(x(t),a',b')_q\u_m^{\alpha',\beta'}(x(t),a',b')_q
\rho'(t)\btu x(t-1/2)=(d_k')^2\delta_{k,m},
\eq
where $\rho'(t)$ and $d_k'$ are the weight function $\rho$ and the
norm $d_n$ given in table \ref{tabla1} with the corresponding change
$a$, $b$, $\alpha$, $\beta$, $s$, $n$ by $a'$, $b'$, $\alpha'$, $\beta'$, $t$, $k$.

Furthermore, with the above choice \refe{change} 
of the parameters of $\u_k^{\alpha',\beta'}(x(t),a',b')_q$,
the hypergeometric function ${}_4 \fhyp_3$ in \refe{dua-rac} coincides with the
function ${}_4 \fhyp_3$ in \refe{pol-rac-nu-F} and therefore the following relation
between the polynomials $\u_k^{\alpha',\beta'}(x(t),a',b')$
and $u_{n}^{\alpha,\beta}(x(s),a,b)_q$ holds
\bq
\label{rel-rac-dua-rac}
\u_k^{\alpha',\beta'}(x(t),a',b')_q= \mathcal{A}(\alpha,\beta,a,b,n,s)
u_{n}^{\alpha,\beta}(x(s),a,b)_q,
\eq
where
$$ \mathcal{A}(\alpha,\beta,a,b,n,s)=
\frac{(-1)^{s-a+n}\tilG(b-a-n)\tilG(s-a+\beta+1)\tilG(b+\alpha+s+1)\tilG(n+1)}
{\tilG(b-s)\tilG(n+\beta+1)\tilG(b+a+\alpha+n+1)\tilG(s-a+1)}.
$$

If we now substitute \refe{rel-rac-dua-rac} in \refe{rel-dual-ort-u} and make the
change \refe{change}, then \refe{rel-dual-ort-u} becomes into relation
\refe{dual-ort-u}, i.e., the polynomial set $\u_k^{\alpha',\beta'}(x(t),a',b')_q$
defined by \refe{dua-rac} (or \refe{rel-rac-dua-rac}) is the dual set associated to
the Racah polynomials $u_{n}^{\alpha,\beta}(x(s),a,b)_q$.

To conclude this study, let us show that the TTRR \refe{ttrr} of the polynomials
$\u_k^{\alpha',\beta'}(x(t),a',b')_q$ is the SODE \refe{eqdif-q-new} of the
polynomials $u_{n}^{\alpha,\beta}(x(s),a,b)_q$ whereas the
SODE \refe{eqdif-q-new} of the  $\u_k^{\alpha',\beta'}(x(t),a',b')_q$
becomes into the  TTRR \refe{ttrr} of $u_{n}^{\alpha,\beta}(x(s),a,b)_q$
and vice versa.

Let denote by $\varsigma(t)$ the $\sigma$ function of the polynomial
$\u_k^{\alpha',\beta'}$, then
$$
\varsigma(t)= [t-a']_{q} [t+b']_{q}[t+a'-\beta']_{q}[b'+\alpha'-t]_{q}=
[n]_q[n+b-a+\alpha+\beta]_q[n+\alpha]_q[b+a-n-\beta]_q,
$$
and therefore,
$$
\varsigma(-t-1)= [\alpha+\beta+n+1]_q[b+a+\alpha+n+1]_q[b-a-n-1]_q[n+\beta+1]_q,
$$
$$
 \lambda_k=[k]_q[\alpha'+\beta'+k+1]_q=[s-a]_q[s+a+1]_q.
$$
For the coefficients $\alpha'_k$, $\beta_k'$ and $\gamma'_k$,
of the TTRR for the polynomials $\u_k^{\alpha',\beta'}$ we
have
$$
\alpha'_k=\dst\frac{[k+1]_q [\alpha'+\beta'+k+1]_q}
{[\alpha'+\beta'+2k+1]_q[\alpha'+\beta'+2k+2]_q}=
\frac{[s-a+1]_q [s+a+1]_q}{[2s+1]_q[2s+2]_q},
$$
$$
\gamma_k'=\frac{[b+\alpha+s]_q[b+\alpha-s]_q[s+a-\beta]_{q}[s-a+\beta]_{q}[b+s]_q[b-s]_q}
{[2s+1]_q[2s]_q},
$$
and
$$
\beta_k'=[n+\frac{\alpha+\beta}2]_q[n+\frac{\alpha+\beta}2+1]_q+
\frac{\sigma(-s-1)}{[2s+1]_q[2s+2]_q}+\frac{\sigma(s)}{[2s+1]_q[2s]_q}.
$$
Also we have $\btu x(t)=[2t+2]_q=[2n+\alpha+\beta+2]_q$ and $x(s)=[s]_q[s+1]_q=
[k+a]_q[k+a+1]_q$.

Let show that the SODE of the Racah polynomials $u_{n}^{\alpha,\beta}(x(s),a,b)_q$
is the TTRR of the $\u_k^{\alpha',\beta'}(x(t),a',b')_q$ polynomials.
First, we substitute the relation \refe{rel-rac-dua-rac} in the SODE
\refe{eqdif-q-new} of the polynomials $u_{n}^{\alpha,\beta}(x(s),a,b)_q$
and use that $u_{n}^{\alpha,\beta}(x(s\pm1),a,b)_q$ is proportional to
$\u_{k\pm1}^{\alpha',\beta'}(x(t),a',b')_q$ (see \refe{rel-rac-dua-rac}).
After some simplification, and using the last formulas we obtain
$$\ba{r}
\alpha'_k \u_{k+1}^{\alpha',\beta'}(x(t),a',b')_q+
\left(\beta_k'-[n]_q[\alpha+\beta+n+1]_q-[\frac{\alpha+\beta}2]_q
[\frac{\alpha+\beta}2+1]_q\right)
\u_k^{\alpha',\beta'}(x(t),a',b')_q\\[3mm]+
\gamma_k'\u_{k-1}^{\alpha',\beta'}(x(t),a',b')_q=0,
\ea
$$
but
$$\ba{c}
[n]_q[\alpha+\beta+n+1]_q+[\frac{\alpha+\beta}2]_q
[\frac{\alpha+\beta}2+1]_q=[n+\frac{\alpha+\beta}2]_q[n+\frac{\alpha+\beta}2+1]_q
=x(t),
\ea
$$
i.e., we obtain the TTRR for the polynomials $\u_k^{\alpha',\beta'}(x(t),a',b')_q$.

If we now substitute \refe{rel-rac-dua-rac} in the TTRR \refe{ttrr}
for the Racah polynomials $u_{n}^{\alpha,\beta}(x(s),a,b)_q$, and
use that $u_{n\pm 1}^{\alpha,\beta}(x(s),a,b)_q  \sim
\u_k^{\alpha',\beta'}(x(t\pm1),a',b')_q$, then we obtain
the SODE
$$\ba{r}\dst
\frac{\varsigma(-t-1)}{\Delta x(t)\Delta x(t-\half)}\u_k^{\alpha',\beta'}(x(t+1),a',b')_q
+\frac{\varsigma(t)}{\nabla x(t) \Delta x(t-\half)}\u_k^{\alpha',\beta'}(x(t-1),a',b')_q
\\[4mm]\dst
-\left[
\frac{\varsigma(-t-1)}{\Delta x(t)\Delta x(t\!-\!\half)}+
\frac{\varsigma(t)}{\nabla x(t) \Delta x(t\!-\!\half)}+[a]_q[a\!+\!1]_q\!-\![k\!+\!a]_q[k\!+\!a\!+\!1]_q
\right]\u_k^{\alpha',\beta'}(x(t),a',b')_q = 0.
\ea
$$
That is  the SODE \refe{eqdif-q-new} of the
$\u_k^{\alpha',\beta'}(x(t),a',b')_q$
since
$$
[a]_q[a+1]_q-[k+a]_q[k+a+1]_q=-[k]_q[k+2a+1]_q=-[k]_q[k+\alpha'+\beta'+1]_q=-\lambda_k.
$$

\subsection{The $q$-Racah polynomials $\widetilde{u}_n^{\alpha,\beta}(x(s),a,b)_q$}
\label{pol-rac-til}
There is another possibility to define the $q$-Racah
polynomials as it is suggested in \cite{mal92,nsu}. It
corresponds to the function
$$
\sigma(s)=[s-a]_{q} [s+b]_{q}[s-a+\beta]_{q}[b+\alpha+s ]_{q},
$$
i.e., $A=1$, $s_1=a$, $s_2=-b$, $s_3=a-\beta$, $s_4=-b-\alpha$.
With this choice we obtain a new family of polynomials
$\widetilde{u}_n^{\alpha,\beta}(x(s),a,b)_q$ that is
orthogonal with respect to the weight function
$$
\mbox{{\normalsize$\rho(s)$}}=
\mbox{\small $ \displaystyle  \frac{\tilG(s+a+1) \tilG(s+a-\beta+1)}
{\tilG(s+\alpha+b+1) \tilG(b+\alpha-s) \tilG(s-a+1)
\tilG(s+b+1)\tilG(s-a+\beta+1) \tilG(b-s)}$} .
$$
All their characteristics can be obtained exactly in the same
way as before. Moreover, they can be also obtained
from the corresponding characteristics of the polynomials
$u_n^{\alpha,\beta}(x(s),a,b)_q$ by changing
$\alpha\to-2b-\alpha$, $\beta\to 2a-\beta $ ---and using the
properties of the functions $\tilG(s)$, $\Gamma_q(s)$,
$(a|q)_n$ and $(a;q)_n$. We will resume the main data of the polynomials
$\widetilde{u}_n^{\alpha,\beta}(x(s),a,b)_q $ in table~\ref{tabla2}.

\begin{table}[t]
\caption{Main data of the {\em q-}Racah polynomials
$\widetilde{u}_n^{\alpha,\beta}(x(s),a,b)_q $\label{tabla2}}

\vspace*{.5cm}

{\small
\begin{center}
{ \renewcommand{\arraystretch}{.35}
\begin{tabular}{|@{}c@{}| | @{}c@{}|}\hline
 &  \\
$P_n(s)$ & $\widetilde{u}_n^{\alpha,\beta}(x(s),a,b)_q   \,,\esp x(s) = [s]_q [s+1]_q$ \\
 &  \\
\hline\hline
 &  \\
$ (a,b)$ & \mbox{
$[a,b-1]$} \\
 &  \\
\hline
 &  \\
$\rho(s)$ & $  \dst \frac{\tilG(s+a+1) \tilG(s+a-\beta+1)}
{\tilG(s+\alpha+b+1) \tilG(b+\alpha-s) \tilG(s-a+1)
\tilG(s+b+1)\tilG(s-a+\beta+1) \tilG(b-s)}  $ \\
 &  \\
 &  \mbox{
$-\half<a\leq b-1, \alpha>-1,-1<\beta<2a+1$} \\
 &  \\
\hline
 &  \\
$\sigma(s)$ & \mbox{
$  [s-a]_{q} [s+b]_{q}[s-a+\beta]_{q}[b+\alpha+s ]_{q} $} \\
 &  \\
\hline
 &  \\
$\sigma(-s-1)$ & \mbox{
$  [s+a+1]_{q} [b-s-1]_{q}[s+a-\beta+1]_{q}[b+\alpha-s-1 ]_{q} $} \\
 &  \\
\hline
 &  \\
$\tau(s) $ & \begin{tabular}{r}
$ [2a-\beta+1]_q [b]_q[b+\alpha]_q-[2b+\alpha-1]_q [a]_q [a-\beta]_q-
[2b+\alpha-1]_q[2a-\beta+1]_q$ \\[3mm]
$-[2b-2a+\alpha+\beta-2]_q x(s) $ \end{tabular}\\
 &  \\
\hline
 &  \\
$\tau_n(s) $ &  \mbox{
\begin{tabular}{r}
$ \!-\![2b\!-\!2a\!+\!\alpha\!+\!\beta\!-\!2n\!-\!2]_q x(s\!+\!\frac n2)\!+\!
[a\!+\!\frac n2\!+\!1]_{q} [b\!-\!\frac n2\!-\!1]_{q}
[a\!+\!\frac n2\!+\!1\!-\!\beta]_{q}[b\!-\!\frac n2\!+\!\alpha\!-\!1]_{q}$
\\[3mm] $-
[a+\frac n2]_{q} [b-\frac n2]_{q}[a+\frac n2-\beta]_{q}[b-\frac n2+\alpha]_{q}$
 \end{tabular}}\\
 &  \\
 \hline
 &  \\
$\lambda _n$ &\mbox{
$[n]_q[2b-2a+\alpha+\beta-n-1]_q$} \\
 &  \\
\hline
 &  \\
$B_n$ &  $  \dst\frac{1}{[n]_q!}$
 \\
 &  \\
 \hline
 &  \\
$d_n^2$&$\dst\frac{\tilG(2a\!+\!n\!-\!\beta\!+\!1)
\tilG(2b\!-\!2a\!+\!\alpha\!+\!\beta\!-\!n)
[2b\!-\!2a\!-\!2n\!-\!1\!+\!\alpha\!+\!\beta]_q^{-1}}
{\tilG(n\!+\!1)\tilG(b\!-\!a\!-\!n)\tilG(b\!-\!a\!-\!n\!+\!\alpha)
\tilG(b\!-\!a\!+\!\beta\!-\!n)\tilG(2b\!+\!\alpha\!-\!n)
\tilG(b\!-\!a\!+\!\alpha\!+\!\beta\!-\!n)}\,\,$\, \\
 &  \\
 \hline
 &  \\
$\,\rho_n(s)\,$ & $ \dst \frac{\tilG(s\!+\!a\!+\!n\!+\!1) \tilG(s\!+\!a\!+\!n\!-\!\beta\!+\!1)}
{\tilG(s\!+\!\alpha\!+\!b\!+\!1) \tilG(b\!+\!\alpha\!-\!s\!-\!n) \tilG(s\!-\!a\!+\!1)
\tilG(s\!+\!b\!+\!1)\tilG(s\!-\!a\!+\!\beta\!+\!1) \tilG(b\!-\!s\!-\!n)}   $  \\
 &  \\
 \hline
  & \\
$a_n$ & $\dst\frac{(-1)^n\tilG[2b-2a+\alpha+\beta-n]_q}
{[n]_q!\tilG[2b-2a+\alpha+\beta-2n]_q}$
\\
 & \\
\hline
 &  \\
$\alpha _n$ &  $  \dst-\frac{ [n+1]_q [2b-2a+\alpha+\beta-n-1]_q}
{[2b-2a+\alpha+\beta-2n-1]_q[2b-2a+\alpha+\beta-2n-2]_q}$ \\
 &  \\
\hline
 &  \\
 $\beta_n $ &  \begin{tabular}{c}
$\dst[a]_q[a\!+\!1]_q \!+\!\frac{[2b\!-\!2a\!+\!\alpha\!+\!\beta\!-\!n\!-\!1]_q
[a\!-\!b\!+\!n\!+\!1]_q[2a\!-\!\beta\!+\!n\!+\!1]_q
[a\!-\!b\!-\!\alpha\!+\!n\!+\!1]_q}
{[2b\!-\!2a\!+\!\alpha\!+\!\beta\!-\!2n\!-\!1]_q[2b\!-\!2a\!+\!\alpha\!+\!\beta\!-\!2n\!-\!2]_q}$ \\[3mm]
$\dst\!+\!\frac{[2b\!+\!\alpha\!-\!n]_q[b\!-\!a\!+\!\alpha\!+\!\beta\!-\!n]_q[b\!-\!a\!+\!\beta\!-\!n]_q[n]_q}
{[2b\!-\!2a\!+\!\alpha\!+\!\beta\!-\!2n\!-\!1]_q[2b\!-\!2a\!+\!\alpha\!+\!\beta\!-\!2n]_q}$
\end{tabular}\\
 &  \\
\hline
 &  \\
$\gamma_n$ & $ -\dst\frac{[2a\!-\!\beta\!+\!n]_q[b\!-\!a\!-\!n]_q
[b\!-\!a\!-\!n\!+\!\alpha\!]_q[b\!-\!a\!-\!n\!+\!\beta\!]_q
[2b\!+\!\alpha\!-\!n]_q[b\!-\!a\!+\!\alpha\!+\!\beta\!-\!n]_q}
{[2b\!-\!2a\!+\!\alpha\!+\!\beta\!-\!2n\!-\!1]_q
[2b\!-\!2a\!+\!\alpha\!+\!\beta\!-\!2n]_q} $ \\
 &  \\
\hline
\end{tabular}  }
\end{center}  }
\end{table}

\subsubsection{The hypergeometric representation}

For the $\widetilde{u}_n^{\alpha,\beta}(x(s),a,b)_q$ polynomials
we have the following hypergeometric representation
\bq
\begin{split}\dst
\widetilde u_n^{\alpha,\beta}(x(s),a,b)_q= & \dst
\frac{q^{-\frac{n}{2}(4a-2b-\alpha-\beta+n+1)}
(q^{a-b+1};q)_n (q^{2a-\beta+1};q)_n (q^{a-b-\alpha+1};q)_n}
{\qdiv^{2n} (q;q)_n}\times\\ &
{}_{4}\serie_3 \left(\ba{c} q^{-n},q^{2a-2b-\alpha-\beta+n+1}, q^{a-s},
q^{a+s+1}  \\ q^{a-b+1},q^{2a-\beta+1},q^{a-b-\alpha+1}\ea
\,\bigg|\, q \,,\, q \right),
\end{split}
\label{pol-til-rac-nu}
\eq
or, in terms of the $q$-hypergeometric series \refe{q-hip-def},
\bq\label{pol-til-rac-nu-F}
\begin{split}\dst
\widetilde u_n^{\alpha,\beta}(x(s),a,b)_q= &
\dst \frac{(a-b+1|q)_n (2a-\beta+1|q)_n
(a-b-\alpha+1|q)_n}{[n]_{q}!}\times\\  &
{}_{4}\fhyp_3 \left(\ba{c} -n,2a-2b-\alpha-\beta+n+1, a-s,a+s+1 \\
a-b+1,2a-\beta+1, a-b-\alpha+1 \ea \,\bigg|\, q \,,\, 1 \right).
\end{split}
\eq
Using the Sears transformation formula \cite[Eq. (III.15)]{gasper1}
we obtain the equivalent representation formulas
\bq
\begin{split}\dst
\widetilde u_n^{\alpha,\beta}(x(s),a,b)_q= & \dst
\frac{q^{-\frac{n}{2}(2a-4b-\alpha-\beta+n+1)}
(q^{a-b+1};q)_n (q^{-2b-\alpha+1};q)_n(q^{-\beta+a-b+1};q)_n}{\qdiv^{2n}
(q;q)_n}\times\\  & {}_{4}\serie_3
\left(\ba{c} q^{-n},q^{2a-2b-\alpha-\beta+n+1}, q^{-b-s},
q^{-b+s+1}  \\ q^{a-b+1},q^{-2b-\alpha+1},q^{a-b-\beta+1}\ea\,\bigg|\, q
\,,\, q \right),
\end{split}
\label{pol-til-rac-nu2}
\eq
and
\bq\label{pol-til-rac-nu-F-2}
\begin{split}\dst
\widetilde u_n^{\alpha,\beta}(x(s),a,b)_q= & \dst \frac{(a-b+1|q)_n (-2b-\alpha+1|q)_n
(a-b-\beta+1|q)_n}{[n]_{q}!}\times\\  &
_{4}\fhyp_3 \left(\ba{c} -n, 2a-2b-\alpha-\beta+n+1, -b-s,-b+s+1 \\
 a-b+1,-2b-\alpha+1, a-b-\beta+1 \ea \,\bigg|\, q \,,\, 1 \right).
\end{split}
\eq
\noindent\textbf{Remark:} From the above formulas follow 
that the polynomials $\widetilde u_n^{\alpha,\beta}(x(s),a,b)_q$
are multiples of the standard $q$-Racah polynomials 
$R_n(\mu(q^{a-s});q^{a-b-\alpha},q^{a-b-\beta},q^{a-b},q^{a+b}|q)$.\\

Moreover, from the above hypergeometric representations follow the values
\bq\label{val-til-a}\begin{split}
\widetilde u_n^{\alpha,\beta}(x(a),a,b)_q=&\dst
\dst \frac{(a-b+1|q)_n (2a-\beta+1|q)_n
(a-b-\alpha+1|q)_n}{[n]_{q}!}\\
=&\frac{(q^{a-b+1};q)_n (q^{2a-\beta+1};q)_n (q^{a-b-\alpha+1};q)_n}
{q^{\frac{n}{2}(4a-2b-\alpha-\beta+n+1)}\qdiv^{2n} (q;q)_n},
\end{split}
\eq
\bq\label{val-til-b}
\begin{split}
\widetilde u_n^{\alpha,\beta}(x(b-1),a,b)_q =&
 \dst \frac{(a-b+1|q)_n (-2b-\alpha+1|q)_n
(a-b-\beta+1|q)_n}{[n]_{q}!}\\  =&\dst
\frac{(q^{a-b+1};q)_n (q^{-2b-\alpha+1};q)_n(q^{-\beta+a-b+1};q)_n}
{q^{\frac{n}{2}(2a-4b-\alpha-\beta+n+1)}\qdiv^{2n}
(q;q)_n}.
\end{split}
\eq

Using   \refe{form-exp-q} we obtain an explicit formula\footnote{Obviously the formulas
(\ref{pol-til-rac-nu}--\ref{pol-til-rac-nu-F-2}) also give two equivalent
explicit formulas.}
$$
\ba{l}
\widetilde u_n^{\alpha,\beta}(x(s),a,b)_q  =  \dst
\frac{ \tilG(s-a+1) \tilG(s+b+1)\tilG(s-a+\beta+1) \tilG(b-s)
\tilG(s+\alpha+b+1)}{\tilG(s+a+1) \tilG(s+a-\beta+1)}\times
\ea
$$
\bq
\ba{l}
\tilG(b+\alpha-s)
\dst \sum_{k=0}^{n} \frac{(-1)^{k+n} [2s+2k-n+1]_q\tilG(s+k+a+1)\tilG(2s+k-n+1)}
{\tilG(k+1) \tilG(n-k+1) \tilG(2s+k+2)\tilG(s-n+k-a+1)\tilG(b-s-k)}\times
\ea\label{q-til-racah}
\eq
$$
\ba{l}
\quad\dst\frac{\tilG(s+k+a-\beta+1)}
{\tilG(s+k-n+\alpha+b+1)\tilG(b+\alpha-s-k)\tilG(s-n+k+b+1)\tilG(s-n+k-a+\beta+1)}.
\ea
$$
From this expression follows that
\bq\ba{l}\dst
\widetilde u_n^{\alpha,\beta}(x(a),a,b)_q =
\frac{ \tilG(b-a)\tilG(2a-\beta+n+1)\tilG(b-a+\alpha)}
{[n]!\tilG(b-a-n)\tilG(2a-\beta+1)\tilG(b-a+\alpha-n)},\\[4mm]
\dst
\widetilde u_n^{\alpha,\beta}(x(b-1),a,b)_q =
\frac{(-1)^n\tilG(b-a)\tilG(2b+\alpha)\tilG(b-a+\beta)}
{[n]!\tilG(b-a-n)\tilG(2b+\alpha-n)\tilG(b-a+\beta-n)},
\ea
\label{val-til-ext}
\eq
that are in agreement with the values \refe{val-til-a} and \refe{val-til-b}
obtained before.

{From} the hypergeometric representation follows the symmetry property
$$
\widetilde u_n^{\alpha,\beta}(x(s),a,b)_q=\widetilde u_n^{-b-a+\beta,b+a+\alpha}(x(s),a,b)_q.
$$

\subsubsection{The  differentiation formulas}

Next we use the differentiation formulas
\refe{for-dif} and \refe{dif-for-2} to obtain
\bq\label{dif-for-til-rac1}
\frac{\btu \widetilde u_n^{\alpha,\beta}(x(s),a,b)_q }{\btu x(s)}=
-[2b-2a+\alpha+\beta-n-1]_q
 \widetilde u_{n-1}^{\alpha,\beta}(x(s+\half),a+\half,b-\half)_q,
\eq
\bq\label{dif-for-til-rac2}\begin{split}
[n]_q[2s+1]_q \widetilde u_n^{\alpha,\beta}(x(s),a,b)_q = &
\sigma(-s-1) \widetilde u_{n-1}^{\alpha,\beta}(x(s+\half),a+\half,b-\half)_q\\[4mm] & -
\sigma(s) \widetilde u_{n-1}^{\alpha,\beta}(x(s-\half),a+\half,b-\half)_q,
\end{split}\eq
respectively. Finally,  the formulas \refe{d3a} (or \refe{d3})
and \refe{d4} lead to the following {\it differentiation formulas}
\bq \label{first-diff-form-til}
\sigma(s)\frac{\nabla \widetilde u_n^{\alpha,\beta}(x(s),a,b)_q }{[2s]_q}=
\!-\!\frac{[2b\!-\!2a\!+\!\alpha\!+\!\beta\!-\!n\!-\!1]_q}
{[2b\!-\!2a\!+\!\alpha\!+\!\beta\!-\!2n\!-\!2]_q}\Bigg[
\tau_n(s)\widetilde u_n^{\alpha,\beta}(x(s),a,b)_q \!-\![n\!+\!1]_q
\widetilde u_{n\!+\!1}^{\alpha,\beta}(x(s),a,b)_q \Bigg],
\eq
\bq \label{scnd-diff-form-til}
\ba{l}\dst
\sigma(-s-1)\frac{\Delta \widetilde{u}_{n}^{\alpha,\beta}(x(s),a,b)_q}{[2s\!+\!2]_q}=
\!-\!\frac{[2b\!-\!2a\!+\!\alpha\!+\!\beta\!-\!n\!-\!1]_q}
{[2b\!-\!2a\!+\!\alpha\!+\!\beta\!-\!2n\!-\!2]_q}\times\\[5mm]
\dst \quad \Bigg[(\tau_n(s)\!+\![n]_q[2b\!-\!2a\!+\!\alpha\!+\!\beta\!-\!2n\!-\!2]_q
[2s\!+\!1]_q)\widetilde u_{n}^{\alpha,\beta}(x(s),a,b)_q \!-\![n\!+\!1]_q
\widetilde u_{n\!+\!1}^{\alpha,\beta}(x(s),a,b)_q
\Bigg],
\ea
\eq
respectively, where  $\tau_n(s)$ is given in  table \ref{tabla2}.

\subsection{The dual set to $\widetilde u_n^{\alpha,\beta}(x(s),a,b)_q$}

To obtain the dual set to $\widetilde u_n^{\alpha,\beta}(x(s),a,b)_q$ we use
the same method as in the previous section. We start from the orthogonality relation
\refe{norm} for the  $\widetilde u_n^{\alpha,\beta}(x(s),a,b)_q$ polynomials
defined by \refe{pol-til-rac-nu-F-2} and write the dual relation
\bq\label{dual-ort-u-til}
\sum_{n=0}^{N-1} \widetilde u_{n}^{\alpha,\beta}(x(s),a,b)_q
\widetilde u_{n}^{\alpha,\beta}(x(s'),a,b)_q
\frac1{d_n^2}=\frac1{\rho(s)\btu x(s-1/2)}\delta_{s,s'},\quad N=b-a,
\eq
where $\rho$ and $d_n^2$ are the weight function and the norm of the
$\widetilde u_{n}^{\alpha,\beta}(x(s),a,b)_q$ given in table
\ref{tabla2}. Furthermore, from \refe{pol-til-rac-nu-F-2} follows that the functions
$\widetilde u_n^{\alpha,\beta}(x(s),a,b)_q$ are polynomials of degree
$k=b-s-1$ on the lattice $x(t)=[t]_q[t+1]_q$ where $t=b-a-n+\frac{\alpha+\beta}2-1$
(the proof is similar to the one presented in section \ref{sec-dual-u} and we will
omit it here). To identify the dual set let us define the new set
\bq\begin{split}
\widetilde\u_k^{\alpha',\beta'}(x(t),a',b')_q=&
\frac{(-1)^{k} \tilG(b'-a')\tilG(b'-a'+\beta')\tilG(2b'+\alpha')}
{[k]!\tilG(b'-a'-k)\tilG(b'-a'+\beta'-k)\tilG(2b'+\alpha'-k)}
\times\\[4mm] &
_{4}\fhyp_3 \left(\ba{c} -k, 2a'-2b'-\alpha'-\beta'+k+1, -b'-t,-b'+t+1 \\
 a'-b'+1,-2b'-\alpha'+1, a'-b'-\beta'+1 \ea \,\bigg|\, q \,,\, 1 \right),
\end{split}
\label{dua-rac-til}
\eq
where
\bq\label{change-til}
k=b-s-1, \quad t=b-a-n+\frac{\alpha+\beta}2-1,\quad
a'=\frac{\alpha+\beta}2, \quad b'=b-a+\frac{\alpha+\beta}2,
\quad \alpha'=2a-\beta, \quad \beta'=\beta.
\eq
Obviously they satisfy the following orthogonality relation
\bq\label{rel-ort-dua-til}
\sum_{t=a'}^{b'-1} \widetilde\u_k^{\alpha',\beta'}(x(t),a',b')_q
\widetilde\u_m^{\alpha',\beta'}(x(t),a',b')_q
\rho'(t)\btu x(t-1/2)=(d_k')^2\delta_{k,m},
\eq
where now $\rho'(t)$ and $d_k'$ are the weight function $\rho$ and the norm
 $d_n$, respectively, given in table \ref{tabla2} with the
corresponding change of the parameters $a,b,\alpha,\beta,n,s$ by
$a',b',\alpha',\beta',k,t$ \refe{change-til}.

Furthermore, with the above definition \refe{change-til}
for the parameters of $\widetilde\u_k^{\alpha',\beta'}(x(t),a',b')_q$,
the hypergeometric function ${}_4 \fhyp_3$ in \refe{dua-rac-til} coincides with the
function ${}_4 \fhyp_3$ in \refe{pol-til-rac-nu-F-2} and therefore the following relation
between the polynomials $\widetilde\u_k^{\alpha',\beta'}(x(t),a',b')$
and $\widetilde u_{n}^{\alpha,\beta}(x(s),a,b)_q$ holds
\bq
\label{rel-rac-dua-rac-til}
\widetilde\u_k^{\alpha',\beta'}(x(t),a',b')_q= \widetilde{\mathcal{A}}(\alpha,\beta,a,b,n,s)
\widetilde u_{n}^{\alpha,\beta}(x(s),a,b)_q,
\eq
where
$$ \widetilde{\mathcal{A}}(\alpha,\beta,a,b,n,s)=
\frac{(-1)^{b-s-1-n}\tilG(b-a-n)\tilG(2b+\alpha-n)\tilG(b-a+\beta-n)\tilG(n+1)}
{\tilG(b-s)\tilG(s-a+\beta+1)\tilG(s+b+\alpha+1)\tilG(s-a+1)}.
$$

To prove that the polynomials $\widetilde\u_k^{\alpha',\beta'}(x(t),a',b')_q$ are
the dual set to $\widetilde u_n^{\alpha,\beta}(x(s),a,b)_q$ it is sufficient to
substitute \refe{rel-rac-dua-rac-til} in \refe{rel-ort-dua-til} and do the change
\refe{change-til} that transforms \refe{rel-ort-dua-til} into \refe{dual-ort-u}.

Let also mention that, as in the case of the $q$-Racah polynomials,
the TTRR \refe{ttrr} of the polynomials  $\widetilde\u_k^{\alpha',\beta'}(x(t),a',b')_q$
is the SODE \refe{eqdif-q-new} of the
polynomials $\widetilde u_{n}^{\alpha,\beta}(x(s),a,b)_q$ whereas the
SODE \refe{eqdif-q-new} of the  $\widetilde\u_k^{\alpha',\beta'}(x(t),a',b')_q$
becomes into the  TTRR \refe{ttrr} of $\widetilde u_{n}^{\alpha,\beta}(x(s),a,b)_q$
and vice versa.

To conclude this section let us point out that there exist a simple relation
connecting both polynomials $u_{n}^{\alpha,\beta}(x(s),a,b)_q$ and
$\tilde u_{n}^{\alpha,\beta}(x(s),a,b)_q$ (see \refe{con-u-util} from below). We will
establish it at the end of the next section.

\section{Connection with the $6j$-symbols of the q-algebra $SU_q(2)$}

\subsection{$6j$-symbols of the quantum algebra $SU_q(2)$}

It is known (see e.g. \cite{smi2} and references therein) that the
Racah coefficients $U_q(j_1\,j_2\,j\,j_3;j_{12}\,j_{23})$ are
used for the transition from the coupling scheme of three angular momenta
$j_1,j_2,j_3$
$$
|j_1j_2(j_{12}),j_3:jm\rangle =\sum_{m_1,m_2,m_3,m_{12}}
\langle j_1m_1j_2m_2|j_{12}m_{12}\rangle \langle j_{12}m_{12}j_3m_3|jm\rangle
|j_1m_1\rangle|j_2m_2\rangle|j_3m_3\rangle,
$$
to the following ones
$$
|j_1j_2j_3(j_{23}):jm\rangle =\sum_{m_1,m_2,m_3,m_{23}}
\langle j_2m_2j_3m_3|j_{23}m_{23}\rangle \langle j_{1}m_{1}j_{23}m_{23}|jm\rangle
|j_1m_1\rangle|j_2m_2\rangle|j_3m_3\rangle,
$$
where $\langle j_am_aj_bm_b|j_{ab}m_{ab}\rangle$ denotes the Clebsh-Gordon
Coefficients of the quantum algebra $su_q(2)$. In fact we have that the recoupling
is given by
$$
|j_1j_2(j_{12}),j_3:jm\rangle=\sum_{j_{23}}U_q(j_1\,j_2\,j\,j_3;j_{12}\,j_{23})
|j_1j_2j_3(j_{23}):jm\rangle.
$$
The Racah coefficients $U$ define an unitary matrix, i.e., they satisfy
the orthogonality relations
\bq
\label{ort-rel1}
\sum_{j_{23}}U_q(j_1\,j_2\,j\,j_3;j_{12}\,j_{23})U_q(j_1\,j_2\,j\,j_3;j_{12}'\,j_{23})
=\delta_{j_{12},j_{12}'},
\eq
\bq
\label{ort-rel2}
\sum_{j_{12}}U_q(j_1\,j_2\,j\,j_3;j_{12}\,j_{23})U_q(j_1\,j_2\,j\,j_3;j_{12}\,j_{23}')
=\delta_{j_{23},j_{23}'}.
\eq
Usually instead of the Racah coefficients is more convenient to use the $6j$-symbols
defined by
$$
U_q(j_1\,j_2\,j\,j_3;j_{12}\,j_{23})=(-1)^{j_1+j_2+j_3+j}\sqrt{[2j_{12}+1]_q[2j_{23}+1]_q}
\left\{\ba{ccc} j_1 & j_2 & j_{12} \\ j_3 & j & j_{23} \ea \right\}_q.
$$

The $6j$-symbols have the following symmetry property
\bq\label{sim-6j}
\left\{\ba{ccc} j_1 & j_2 & j_{12} \\ j_3 & j & j_{23} \ea \right\}_q=
\left\{\ba{ccc} j_3 & j_2 & j_{23} \\ j_1 & j & j_{12} \ea \right\}_q.
\eq

Here and without lost of generality we will suppose that $j_1\geq j_2$ and $j_3\geq j_2$,
then for the moments $j_{23}$ and $j_{12}$ we have the intervals
$$
j_3-j_2\leq j_{23}\leq j_2+j_3,\quad j_1-j_2\leq j_{12}\leq j_1+j_2,
$$
respectively. Now, in order to avoid any other restrictions on these two
momenta (caused by the so called triangle inequalities for the $6j$-symbols)
we will assume that the following restrictions hold
$$
|j-j_3|\leq \min(j_{12})=j_1-j_2,
\quad
|j-j_1|\leq \min(j_{23})=j_3-j_2.
$$

\subsection{$6j$-symbols and the $q$-Racah polynomials
$u_n^{\alpha,\beta}(x(s),a,b)_q$}

Now we are ready to establish the connection of $6j$-symbols with the $q$-Racah
polynomials. We fix the variable $s$  as $s=j_{23}$ that runs on the interval
$a\leq s\leq b-1$ where $a=j_3-j_2$, $b=j_2+j_3+1$. Let us put
\bq
\label{6j-rac}
(-1)^{j_1+j_{23}+j}\sqrt{[2j_{12}+1]_q}
\left\{\ba{ccc} j_1 & j_2 & j_{12} \\ j_3 & j & j_{23} \ea \right\}_q=
\sqrt{\frac{\rho(s)}{d_n^2}} \, u_{n}^{\alpha,\beta}(x(s),a,b)_q,
\eq
where $\rho(s)$ and $d_n$ are the weight function and the norm, respectively, of the
$q$-Racah polynomials on the lattice \refe{red-cua} $u_{n}^{\alpha,\beta}(x(s),a,b)_q$,
and $n=j_{12}-j_1+j_2$, $\alpha=j_1-j_2-j_3+j\geq0$,
$\beta=j_1-j_2+j_3-j\geq0$.\footnote{Notice that this is equivalent to the following
setting
\[\begin{split}
&j_1=(b-a-1+ \alpha + \beta)/2
,\quad
j_2=(b-a-1)/2
,\quad
j_3=(a+b-1)/2
,\\
&j_{12}=(2 n + \alpha + \beta)/2
,\quad
j_{23}= s,\quad
j= (a+ b-1+\alpha - \beta)/2.
\end{split}
\]}

To verify the above relation we use the recurrence relation \cite[Eq. (5.17)]{smi3} \bq\label{rr-6j}\begin{split} [2]_q& [2j_{23} +2]_q
A_q^- \left\{\ba{ccc} j_1 & j_2 & j_{12} \\ j_3 & j & j_{23}-1 \ea
\right\}_q\!\!\!-\\ &
\Big(([2j_{23}]_q[2j_1+2]_q\!-\![2]_q[j-j_{23}+j_1+1]_q[j+j_{23}-j_1]_q)\times
\\ &
([2j_2]_q[2j_{23}+2]_q-[2]_q[j_3-j_2+j_{23}+1]_q[j_3+j_2-j_{23}]_q)-\\
&
([2j_2]_q[2j_{1}+2]_q-[2]_q[j_{12}-j_2+j_1+1]_q[j_{12}+j_2-j_1]_q)[2j_{23}+2]_q
[2j_{23}]_q       
\Big)\times \\  [2j_{23}+1]_q& \left\{\ba{ccc} j_1 & j_2 & j_{12}
\\ j_3 & j & j_{23} \ea \right\}_q+ [2]_q[2j_{23}]_q
A_q^+\left\{\ba{ccc} j_1 & j_2 & j_{12} \\ j_3 & j & j_{23}+1 \ea
\right\}_q=0,
\end{split}
\eq
where
\bq\label{A}
\begin{split}
A_q^-=&
\sqrt{[j+j_{23}+j_1+1]_q[j+j_{23}-j_1]_q[j-j_{23}+j_1+1]_q[j_{23}-j+j_1]_q}\times\\
&\sqrt{[j_2+j_3+j_{23}+1]_q[j_2+j_3-j_{23}+1]_q[j_3-j_2+j_{23}]_q[j_2-j_3+j_{23}]_q},\\
A_q^+=&
\sqrt{[j+j_{23}+j_1+2]_q[j+j_{23}-j_1+1]_q[j-j_{23}+j_1]_q[j_{23}-j+j_1+1]_q}\times\\
&\sqrt{[j_2+j_3+j_{23}+2]_q[j_2+j_3-j_{23}]_q[j_3-j_2+j_{23}+1]_q[j_2-j_3+j_{23}+1]_q}.
\end{split}
\eq
Notice that
$$
A_q^-=\sqrt{\sigma(j_{23})\sigma(-j_{23})},\quad
A_q^+=\sqrt{\sigma(j_{23}+1)\sigma(-j_{23}-1)},
$$
where
\[\begin{split}
\sigma(j_{23})&= [j_{23}-j_3+j_2]_q[j_{23}+j_2+j_3+1]_q[j_{23}-j_1+j]_q[j+j_1-j_{23}+1]_q,\\
\sigma(-j_{23}-1)&= [j_{23}+j_3-j_2+1]_q[j_2+j_3-j_{23}]_q[j_{23}+j_1-j+1]_q[j+j_1+j_{23}+2]_q.
\end{split}
\]

Substituting \refe{6j-rac} in \refe{rr-6j} and simplifying the obtained expression
we get
\[\begin{split}
[2&s]_q\sigma(-s-1)u_{n}^{\alpha,\beta}(x(s+1),a,b)_q+
[2s+2]_q\sigma(s)u_{n}^{\alpha,\beta}(x(s-1),a,b)_q+\\
&\Big(\lambda_n [2s]_q[2s+1]_q[2s+2]_q-[2s]_q\sigma(-s-1)-[2s+2]_q\sigma(s)
\Big)u_{n}^{\alpha,\beta}(x(s),a,b)_q=0,
\end{split}
\]
which is the difference equation for the $q$-Racah polynomials
\refe{eqdif-q-new}. Since $u_{0}^{\alpha,\beta}(x(s),a,b)_q=1$, \refe{6j-rac}
leads to
$$
(-1)^{j_1+j_{23}+j}\sqrt{[2j_{1}-2j_2+1]_q}
\left\{\ba{ccc} j_1 & j_2 & j_{1}-j_2 \\ j_3 & j & j_{23} \ea \right\}_q=
\sqrt{\frac{\rho(s)}{d_0^2}}\quad \Rightarrow
$$
\[
\begin{split}&
\left\{\ba{ccc} j_1 & j_2 & j_{1}-j_2 \\ j_3 & j & j_{23} \ea \right\}_q:=
\left\{\ba{ccc} j_1 & j_2 & j_{1}-j_2 \\ j_3 & j & s \ea \right\}_q\\
=& (-1)^{j+j_1+s}\sqrt{\frac{[j_1+j+s+1]_q![j_1+j-s]_q![j_1-j+s]_q![j_3-j_2+s]_q!}
{[j-j_1+s]_q![j_3+j_2-s]_q![j_2-j_3+s]_q![j_2+j_3+s+1]_q!}}\times\\
&\sqrt{\frac{[2j_1-2j_2]_q![2j_2]_q![j_2+j_3+j-j_1]_q!}
{[2j_1+1]_q![j_1+j_3-j_2-j]_q![j_1-j_3-j_2+j]_q![j_1+j_3-j_2+j+1]_q!}}.
\end{split}
\]
Furthermore, substituting the values $s=a$ and $s=b-1$ in \refe{6j-rac} and using
\refe{val-ext} we find
\bq \label{6j-en_a}
\begin{split} & \left\{\ba{ccc} j_1 & j_2 & j_{12} \\ j_3 & j & j_3-j_2 \ea\right\}_q
=(-1)^{j_{12}+j_{3}+j} \times
\\  & \qquad\sqrt{\frac{\qnf{j_{12}+j_3-j} \qnf{2j_2} \qnf{j_{12}+j_3+j+1}
\qnf{2j_3-2j_2}\qnf{j_2-j_1+j_{12}}}
{\qnf{j_1-j_2+j_3-j} \qnf{j_1+j_2-j_{12}} \qnf{j_1-j_2+j_3+j+1}
}} \times \\
& \qquad\sqrt{\frac{\qnf{j_1+j_2-j_3+j} \qnf{j_1-j_2+j_{12}} \qnf{j_3-j_{12}+j}}
{\qnf{2j_3+1} \qnf{j_3-j_1-j_2+j} \qnf{j_{12}-j_3+j} \qnf{j_1+j_2+j_{12}+1}}
}
\end{split}
\eq
and
\bq \label{6j-en_b-1}
\begin{split} & \left\{\ba{ccc} j_1 & j_2 & j_{12} \\ j_3 & j & j_2+j_3 \ea \right\}_q=
(-1)^{j_1+j_2+j_3+j} \times
\\  & \qquad \sqrt{\frac{\qnf{2j_2} \qnf{j_{12}-j_3+j} \qnf{j_2-j_1+j_3+j} \qnf{2j_3} \qnf{j_1+j_2+j_3-j}}
{\qnf{j_1+j_2-j_{12}} \qnf{j_1-j_2-j_3+j} \qnf{j_3-j_{12}+j}
}} \times \\
& \qquad\sqrt{\frac{\qnf{j_2-j_1+j_{12}}\qnf{j_1-j_2+j_{12}} \qnf{j_1+j_2+j_3+j+1}}
{\qnf{2j_2+2j_3+1} \qnf{j_{12}+j_3-j} \qnf{j_1+j_2+j_{12}+1} \qnf{j_{12}+j_3+j+1}}
},
\end{split}
\eq
that are in agreement with the results in \cite{smi2}.

The relation \refe{6j-rac} allows us to obtain  several
recurrence relations for the $6j$-symbols of the quantum algebra $SU_q(2)$
by using the  properties of the $q$-Racah polynomials. So, the TTRR \refe{ttrr} gives
\bq
\label{ttrr-6j-j12}\begin{split}
&\quad [2 j_{12}]_q   \, \widetilde{A}_q^+
\left\{\ba{ccc} j_1 & j_2 & j_{12}+1 \\ j_3 & j & j_{23} \ea \right\}_q+[2 j_{12}+2]_q  \, \widetilde{A}_q^-
\left\{\ba{ccc} j_1 & j_2 & j_{12}-1 \\ j_3 & j & j_{23} \ea \right\}_q\\
&-\Big( [2j_{12}]_q [2j_{12}+1]_q [2j_{12}+2]_q \big([j_{23}]_q
[j_{23}+1]_q-[j_3-j_2]_q [j_3-j_2+1]_q\big)+[2j_{12}]_q
\times \\ & [j_1-j_2+j_{12}+1]_q [j_{12}-j_1-j_2]_q [j_{12}+j_3-j+1]_q
[j_{12}+j_3+j+2]_q -[2j_{12}+2]_q \times
\\ & [j_{12}-j_3+j]_q [j_1+j_2+j_{12}+1]_q [j_3-j_{12}+j+1]_q
[j_2-j_1+j_{12}]_q\Big)\left\{\ba{ccc} j_1 & j_2 & j_{12} \\ j_3 & j & j_{23} \ea \right\}_q =0,
\end{split}
\eq
where
\bq
\label{A-til}
\ba{rl}
\widetilde{A}_q^- =&  \dst \sqrt{[j_2-j_1+j_{12}]_q [j_1-j_2+j_{12}]_q [j_{12}-j_3+j]_q [j_{12}+j_3-j]_q
[j_1+j_2+j_{12}+1]_q } \times  \\ & \dst \sqrt{[j_{12}+j_3+j+1]_q [j_1+j_2-j_{12}+1]_q [j_3-j_{12}+j+1]_q}
\\
\widetilde{A}_q^+ =&  \dst \sqrt{[j_2-j_1+j_{12}+1]_q [j_1-j_2+j_{12}+1]_q
[j_{12}-j_3+j+1]_q [j_{12}+j_3-j+1]_q}  \times  \\ &
\dst \sqrt{[j_1+j_2+j_{12}+2]_q [j_{12}+j_3+j+2]_q [j_1+j_2-j_{12}]_q [j_3-j_{12}+j]_q}.
\ea
\eq
The expressions \refe{dif-for-rac1} and \refe{dif-for-rac2} yield
\bq \label{6j-for-dif}
\begin{split}
& \sqrt{\sigma(j_{23}+1)}
\left\{\ba{ccc} j_1 & j_2 & j_{12} \\ j_3 & j & j_{23}+1 \ea \right\}_q 
+\sqrt{ \sigma(-j_{23}-1)}
\left\{\ba{ccc} j_1 & j_2 & j_{12} \\ j_3 & j & j_{23}\ea \right\}_q
\\ & =[2j_{23}+2]_q \sqrt{[j_2-j_1+j_{12}]_q[j_1-j_2+j_{12}+1]_q}
\left\{\ba{ccc} j_1+\half & j_2-\half & j_{12} \\ j_3 & j & j_{23}+\half \ea \right\}_q, \end{split}
\eq
and
\bq \label{6j-dif-for-2}
\begin{split}
\sqrt{\sigma(-j_{23}-1)}&
\left\{\ba{ccc} j_1+\half & j_2-\half & j_{12} \\ j_3 & j & j_{23}+\half \ea \right\}_q  +
\sqrt{\sigma(j_{23})}\left\{\ba{ccc} j_1+\half & j_2-\half & j_{12} \\ j_3 & j & j_{23}-\half\ea \right\}_q
\\ & = [2j_{23}+1]_q \sqrt{[j_{12}-j_1+j_2]_q[j_{12}+j_1-j_2+1]_q}
\left\{\ba{ccc} j_1 & j_2 & j_{12} \\ j_3 & j & j_{23}\ea \right\}_q,
\end{split}
\eq
respectively, whereas the differentiation formulas \refe{first-diff-form}--\refe{scnd-diff-form}
give
\bq
\label{rr1-dif-6j}
\begin{split}
&[2j_{12}+2]_q A_q^-\left\{\ba{ccc} j_1 & j_2 & j_{12} \\ j_3 & j & j_{23}-1 \ea \right\}_q+
[2j_{23}]_q \widetilde{A}_q^+ \left\{\ba{ccc} j_1 & j_2 & j_{12}+1 \\ j_3 & j & j_{23} \ea \right\}_q+
\\
&\Big(\sigma(j_{23})[2j_{12}+2]_q+ [j_1-j_2+j_{12}+1]_q [2j_{23}]_q \Lambda(j_{12},j_{23},j_1,j_2)
\left\{\ba{ccc} j_1 & j_2 & j_{12} \\ j_3 & j & j_{23} \ea \right\}_q=0\\
&\end{split}
\eq
and
\bq
\label{rr2-dif-6j}
\begin{split}
&[2j_{12}+2]_q  A_q^+\left\{\ba{ccc} j_1 & j_2 & j_{12} \\ j_3 & j & j_{23}+1 \ea \right\}_q
-[2j_{23}+2]_q \widetilde{A}_q^+
\left\{\ba{ccc} j_1 & j_2 & j_{12}+1 \\ j_3 & j & j_{23} \ea \right\}_q+\\
&\Big([2j_{12}+2]_q \sigma(-j_{23}-1)-[2j_{23}+2]_q[j_1-j_2+j_{12}+1]_q
\left( \Lambda(j_{12},j_{23},j_1,j_2)+ \right.\\
& \quad \left.[j_{12}-j_1+j_2]_q[2j_{12}+2]_q[2j_{23}+1]_q\right)
\left\{\ba{ccc} j_1 & j_2 & j_{12} \\ j_3 & j & j_{23} \ea \right\}_q=0,
\end{split}
\eq
respectively, where $A_q^{\pm}$ are given by \refe{A}, $\widetilde{A}_q^{\pm}$ by
\refe{A-til} and
\[\begin{split}
\Lambda(j_{12},j_{23},j_1,j_2)=& \sigma \left(\frac{-j_{12}+j_1-j_2}2-1\right)-
\sigma\left(\frac{-j_{12}+j_1-j_2}2\right)- \\
&[2j_{12}+2]_q\left[j_{23}+\frac{j_{12}-j_1+j_2}2\right]_q\left[j_{23}+
\frac{j_{12}-j_1+j_2}2+1\right]_q.
\end{split}
\]

Using the hypergeometric representations \refe{pol-rac-nu-F} and \refe{pol-rac-nu-F-2}
we obtain the representation of the $6j$-symbols in terms of the $q$-hypergeometric
function\footnote{To obtain the representation in terms of the basic hypergeometric
series it is sufficient to use the relation \refe{ser-q-fhyp}.} \refe{q-hip-def}
\[\begin{split}
& \left\{\ba{ccc} j_1 & j_2 & j_{12} \\ j_3 & j & j_{23} \ea \right\}_q
=(-1)^{j_{12}+j_{23}+j_2+j}\frac{\qnf{2j_2}}{\qnf{j_1-j_2+j_3-j}\qnf{j_1-j_2+j_3+j+1}}\\
&\qquad \sqrt{\frac{[j_1+j+j_{23}+1]_q![j_1+j-j_{23}]_q![j_1-j+j_{23}]_q![j_3-j_2+j_{23}]_q!}
{[j-j_1+j_{23}]_q![j_3+j_2-j_{23}]_q![j_2-j_3+j_{23}]_q![j_2+j_3+j_{23}+1]_q!}}\times\\
&\qquad \sqrt{\frac{[j_{12}-j_1+j_2]_q! [j_{12}+j_1-j_2]_q![j_3+j-j_{12}]_q!
[j_3+j_{12}-j]_q![j_3+j_{12}+j+1]_q!}
{[j_{12}-j_3+j]_q! [j_1+j_2+j_{12}+1]_q! [j_1+j_2-j_ {12}]_q!}}\times\\
& \quad\qquad
{}_4\fhyp_3 \left(\ba{c} j_1-j_2-j_{12}, j_1-j_2+j_{12}+1,j_3-j_2-j_{23},j_{23}+j_3-j_2+1 \\
-2j_2,\,\, j_1-j_2+j_3-j+1,\,\, j_1-j_2+j_3+j+2\ea \,\bigg|\, q \,,\, 1 \right),
\end{split}
\]
and
\[\begin{split}
& \left\{\ba{ccc} j_1 & j_2 & j_{12} \\ j_3 & j & j_{23} \ea \right\}_q
=(-1)^{j_{1}+j_{23}+j}\frac{\qnf{2j_2}\qnf{j_2+j_3-j_1+j}}{\qnf{j_1-j_2-j_3+j}}\\
&\qquad \sqrt{\frac{[j_1+j+j_{23}+1]_q![j_1+j-j_{23}]_q![j_1-j+j_{23}]_q![j_3-j_2+j_{23}]_q!}
{[j-j_1+j_{23}]_q![j_3+j_2-j_{23}]_q![j_2-j_3+j_{23}]_q![j_2+j_3+j_{23}+1]_q!}}\times\\
&\qquad \sqrt{\frac{[j_{12}-j_1+j_2]_q! [j_{12}+j_1-j_2]_q!  [j_{12}-j_3+j]_q!}
{[j_1+j_2+j_{12}+1]_q! [j_1+j_2-j_ {12}]_q!  [j_3+j-j_{12}]_q!
[j_{12}+j_3-j]_q![j_3+j_{12}+j+1]_q!}
}\times\\
& \quad\qquad
{}_4\fhyp_3 \left(\ba{c} j_1-j_2-j_{12}, j_1-j_2+j_{12}+1,-j_3-j_2+j_{23},
-j_{23}-j_3-j_2-1 \\
-2j_2,\,\, j_1-j_2-j_3+j+1,\,\, j_1-j_2-j_3-j\ea \,\bigg|\, q \,,\, 1 \right).
\end{split}
\]
Notice that from the above representations the values
\refe{6j-en_a} and \refe{6j-en_b-1} immediately follows.
Notice also that the above formulas give two alternative
explicit formulas for computing the $6j$-symbols. A third explicit formula
follows from  \refe{exp-for-rac}
\[\begin{split}&\hspace{-.75cm}
\left\{\ba{ccc} j_1 & j_2 & j_{12} \\ j_3 & j & j_{23} \ea \right\}_q=
\sqrt{\frac{\qnf{j_{23}+j_2-j_3} \qnf{j_{23}+j_2+j_3+1} \qnf{j_{23}+j-j_1} \qnf{j_2+j_3-j_{23}}}
{\qnf{j_{23}+j_3-j_2} \qnf{j_{23}+j_1-j} \qnf{j_{23}+j_1+j+1} \qnf{j_1+j-j_{23}}}} \times \\
\qquad  &\qquad \qquad \qquad
\sqrt{\frac{\qnf{j_{12}-j_1+j_2} \qnf{j_1-j_2+j_{12}} \qnf{j_1+j_2-j_{12}} \qnf{j_3-j_{12}+j}}
{\qnf{j_{12}-j_3+j} \qnf{j_{12}+j_3-j} \qnf{j_1+j_2+j_{12}+1} \qnf{j_{12}+j_3+j}}} \times \\ &
\sum_{k=0}^{j_{12}-j_1+j_2} \frac{(-1)^{k+j_1+j_{23}+j} [2k+j_1-j_2-j_{12}+2j_{23}+1]_q \qnf{k+j_{23}+j_3-j_2}}
{\qnf{k} \qnf{j_{12}-j_1+j_2-k} \qnf{2j_3+1+k} \qnf{k+j_{23}+j_1-j_{12}-j_3}} \times
\\ & \frac{\qnf{2j_{23}+k-j_{12}+j_1-j_2} \qnf{k+j_{23}+j_1-j} \qnf{k+j_{23}+j_1+j+1} \qnf{j_1+j-j_{23}-k}}
{\qnf{k+j_{23}+j_1-j_{12}+j_3+1} \qnf{k+j_{23}+j-j_2-j_{12}} \qnf{j_2+j_3-j_{23}+1-k}}
\end{split}
\]

To conclude this section let us point out that the orthogonality relations \refe{ort-rel1}
and \refe{ort-rel2} lead to the orthogonality relations for the Racah polynomials
$u_n^{\alpha,\beta}(x(s),a,b)_q$ \refe{pol-rac-nu-F} and their {\em duals}
$\u_k^{\alpha',\beta'}(x(t),a',b')_q$, respectively, and also that the relation
\refe{rel-rac-dua-rac} between $q$-Racah and dual $q$-Racah corresponds to
the symmetry property \refe{sim-6j}.


\subsection{$6j$-symbols and the alternative $q$-Racah polynomials
$\widetilde u_n^{\alpha,\beta}(x(s),a,b)_q$}

In this section we will provide the same comparative analysis but for the
alternative $q$-Racah polynomials $\widetilde u_n^{\alpha,\beta}(x(s),a,b)_q$.
We again choose $s=j_{23}$ that runs on the interval
$[a,b-1]$, $a=j_3-j_2$, $b=j_2+j_3+1$. In this case the connection
is given by formula
\bq
\label{6j-rac-til}
(-1)^{j_{12}+j_3+j} \sqrt{[2j_{12}+1]_q}
\left\{\ba{ccc} j_1 & j_2 & j_{12} \\ j_3 & j & j_{23} \ea \right\}_q=
\sqrt{\frac{\rho(s)}{d_n^2}} \, \widetilde
u_{n}^{\alpha,\beta}(x(s),a,b)_q,
\eq
where $\rho(s)$ and $d_n$ are the weight function and the norm, respectively, of the
alternative $q$-Racah polynomials $\widetilde u_{n}^{\alpha,\beta}(x(s),a,b)_q$
(see Section \ref{pol-rac-til}) on the lattice \refe{red-cua}, and $n=j_1+j_2-j_{12}$,
$\alpha=j_1 - j_2 - j_3 +j \geq0$, $\beta=j_1 - j_2 + j_3- j \geq0$.

Using the above relations we see that the SODE \refe{eqdif-q-new} for the
$\widetilde u_{n}^{\alpha,\beta}(x(s),a,b)_q$ polynomials becomes into the
recurrence relation \refe{rr-6j} as well as the TTRR \refe{ttrr}
becomes into the recurrence relation \refe{ttrr-6j-j12}. Evaluating
\refe{6j-rac-til} in $s=j_{23}=j_3-j_2$ and $s=j_{23}=j_2+j_3+1$  and using \refe{val-til-ext}
we recover the values \refe{6j-en_a} and \refe{6j-en_b-1}, respectively. If we now put
$n=0$, i.e., $j_{12}=j_1+j_2$ we obtain the value
\[
\begin{split}&
\left\{\ba{ccc} j_1 & j_2 & j_{1}+j_2 \\ j_3 & j & j_{23} \ea \right\}_q:=
\left\{\ba{ccc} j_1 & j_2 & j_{1}+j_2 \\ j_3 & j & s \ea \right\}_q\\
& =(-1)^{j_1+j_{2}+j_3+j}\sqrt{\frac{[2j_1]_q![2j_2]_q![j_1+j_2+j_3+j+1]_q![j_1+j_2-j_3+j]_q!}
{[2j_1+2j_2+1]_q![-j_1-j_2+j_3+j]_q![j_2+j_3+s+1]_q!}}\times\\
&\sqrt{\frac{[s-j_1+j]_q![s-j_2+j_3]_q!}
{[j_1+j-s]_q![j_1-j+s]_q![j_1+j+s+1]_q![j_2+j_3-s]_q![j_2-j_3+s]_q!}}.
\end{split}
\]

The expressions \refe{dif-for-til-rac1} and \refe{dif-for-til-rac2} yield
\bq \label{6j-for-dif-til}
\begin{split}
& \sqrt{\varsigma(j_{23}+1)}
\left\{\ba{ccc} j_1 & j_2 & j_{12} \\ j_3 & j & j_{23}+1 \ea \right\}_q 
-\sqrt{ \varsigma(-j_{23}-1)}
\left\{\ba{ccc} j_1 & j_2 & j_{12} \\ j_3 & j & j_{23}\ea \right\}_q
\\ & =[2j_{23}+2]_q \sqrt{[j_1+j_2-j_{12}]_q[j_1+j_2+j_{12}+1]_q}
\left\{\ba{ccc} j_1-\half & j_2-\half & j_{12} \\ j_3 & j & j_{23}+\half \ea \right\}_q,
\end{split}
\eq
and
\bq \label{6j-dif-for-2-til}
\begin{split}
\sqrt{\varsigma(-j_{23}-1)}&
\left\{\ba{ccc} j_1-\half & j_2-\half & j_{12} \\ j_3 & j & j_{23}+\half \ea \right\}_q  -
\sqrt{\varsigma(j_{23})}\left\{\ba{ccc} j_1-\half & j_2-\half & j_{12} \\
j_3 & j & j_{23}-\half\ea \right\}_q
\\ & = [2j_{23}+1]_q \sqrt{[j_{1}+j_2-j_{12}]_q[j_1+j_2+j_{12}+1]_q}
\left\{\ba{ccc} j_1 & j_2 & j_{12} \\ j_3 & j & j_{23}\ea \right\}_q,
\end{split}
\eq
respectively, where
\[\begin{split}
\varsigma(j_{23})&= [j_{23}-j_3+j_2]_q[j_{23}+j_2+j_3+1]_q  [j_{23}-j_1+j+1]_q[j+j_1+j_{23}+1]_q,\\
\varsigma(-j_{23}-1)&= [j_{23}+j_3-j_2+1]_q[j_2+j_3-j_{23}]_q  [j_{23}+j_1-j+1]_q[j+j_1-j_{23}]_q.
\end{split}
\]

The differentiation formulas \refe{first-diff-form-til}--\refe{scnd-diff-form-til} give
\bq
\label{rr1-dif-6j-til}
\begin{split}
&[2j_{12}]_q A_q^-\left\{\ba{ccc} j_1 & j_2 & j_{12} \\ j_3 & j & j_{23}-1 \ea \right\}_q-
[2j_{23}]_q \widetilde{A}_q^- \left\{\ba{ccc} j_1 & j_2 & j_{12}-1 \\ j_3 & j & j_{23} \ea \right\}_q-
\\
&\Big(\varsigma(j_{23})[2j_{12}]_q+ [j_1+j_2+j_{12}+1]_q [2j_{23}]_q \widetilde\Lambda(j_{12},j_{23},j_1,j_2)\Big)
\left\{\ba{ccc} j_1 & j_2 & j_{12} \\ j_3 & j & j_{23} \ea \right\}_q=0\\
&\end{split}
\eq
and
\bq
\label{rr2-dif-6j-til}
\begin{split}
&[2j_{12}]_q  A_q^+\left\{\ba{ccc} j_1 & j_2 & j_{12} \\ j_3 & j & j_{23}+1 \ea \right\}_q
+[2j_{23}+2]_q \widetilde{A}_q^-
\left\{\ba{ccc} j_1 & j_2 & j_{12}-1 \\ j_3 & j & j_{23} \ea \right\}_q-\\
&\Big([2j_{12}]_q \varsigma(-j_{23}-1)-[2j_{23}+2]_q[j_1+j_2+j_{12}+1]_q
\left( \widetilde\Lambda(j_{12},j_{23},j_1,j_2)+ \right.\\
& \quad \left.[j_{1}+j_2-j_{12}]_q[2j_{12}]_q[2j_{23}+1]_q\right)\Big)\Big)
\left\{\ba{ccc} j_1 & j_2 & j_{12} \\ j_3 & j & j_{23} \ea \right\}_q=0,
\end{split}
\eq
respectively, where $A_q^{\pm}$ are given by \refe{A}, $\widetilde{A}_q^{\pm}$ by
\refe{A-til} and
\[\begin{split}
\widetilde\Lambda(j_{12},j_{23},j_1,j_2)=& \varsigma \left(\frac{j_{12}-j_1-j_2}2-1\right)-
\varsigma\left(\frac{j_{12}-j_1-j_2}2\right)- \\
&[2j_{12}]_q\left[j_{23}+\frac{j_1+j_2-j_{12}}2\right]_q\left[j_{23}+
\frac{j_{1}+j_2-j_{12}}2+1\right]_q.
\end{split}
\]

If we now use the hypergeometric representations \refe{pol-til-rac-nu-F} and \refe{pol-til-rac-nu-F-2}
we obtain two new representations of the $6j$-symbols in terms of the $q$-hypergeometric function
\refe{q-hip-def}
\[\begin{split}
& \left\{\ba{ccc} j_1 & j_2 & j_{12} \\ j_3 & j & j_{23} \ea \right\}_q
=(-1)^{j_{12}+j_{3}+j}\frac{\qnf{2j_2}\qnf{j_1+j_2-j_3+j}}{\qnf{j_3-j_2-j_1+j}}\times
\\
&\qquad \sqrt{\frac{[j-j_1+j_{23}]_q![j_3-j_2+j_{23}]_q!}
{[j_1+j+j_{23}+1]_q![j_1+j-j_{23}]_q![j_2-j_3+j_{23}]_q![j_2+j_3+j_{23}+1]_q!
[j_1-j+j_{23}]_q!}}\times\\
&\qquad \sqrt{\frac{[j_3-j_{12}+j]_q! [j_{12}+j_3-j]_q![j_3+j_{12}+j+1]_q![j_1-j_2+j_ {12}+1]_q!}
{[j_3+j_2-j_{23}]_q![j_1+j_2-j_{12}]_q![j_2-j_1+j_{12}]_q![j_{12}-j_3+j]_q![j_1+j_2+j_{12}+1]_q!}}\times\\
& \quad\qquad
{}_4\fhyp_3 \left(\ba{c} j_{12}-j_1-j_2, -j_1-j_2-j_{12}-1,j_3-j_2-j_{23},j_{23}+j_3-j_2+1 \\
-2j_2,\,\, j_3-j_1-j_2+j+1,\,\, j_3-j_1-j_2-j\ea \,\bigg|\, q \,,\, 1 \right),
\end{split}
\]
and
\[\begin{split}
& \left\{\ba{ccc} j_1 & j_2 & j_{12} \\ j_3 & j & j_{23} \ea \right\}_q
=(-1)^{j_{1}+j_{23}+j}
\frac{\qnf{2j_2}\qnf{j_1+j_2+j_3+j}\qnf{j_1+j_2+j_3-j}}{\sqrt{[j_3+j_2-j_{23}]_q! [j_1+j_2-j_{12}]_q!}}\times
\\
&\qquad \sqrt{\frac{[j-j_1+j_{23}]_q![j_3-j_2+j_{23}]_q!}
{[j_1+j+j_{23}+1]_q![j_1+j-j_{23}]_q![j_2-j_3+j_{23}]_q![j_2+j_3+j_{23}+1]_q!
[j_1-j+j_{23}]_q!}}\times\\
&\qquad \sqrt{\frac{[j_{12}-j_3+j]_q! [j_{12}+j_1-j_2+1]_q!}
{[j_3-j_{12}+j]_q![j_{12}+j_3-j]_q![j_1+j_2+j_{12}+1]_q! [j_2-j_1+j_ {12}]_q![j_3+j_{12}+j+1]_q!}}\times\\
& \quad\qquad
{}_4\fhyp_3 \left(\ba{c} j_{12}-j_1-j_2, -j_1-j_2-j_{12}-1,-j_3-j_2-j_{23}-1,
j_{23}-j_3-j_2 \\
-2j_2,\,\, -j_1-j_2-j_3-j,\,\, j-j_1-j_2-j_3\ea \,\bigg|\, q \,,\, 1 \right).
\end{split}
\]
Notice that from the above representations the values
\refe{6j-en_a} and \refe{6j-en_b-1} also follows. Obviously
the above formulas give another two alternative
explicit formulas for computing the $6j$-symbols. Finally,
from  \refe{q-til-racah}
\begin{small}
\[
\begin{split}
& \left\{\ba{ccc} j_1 & j_2 & j_{12} \\ j_3 & j & j_{23} \ea \right\}_q=
\sqrt{\frac{\qnf{j_{23}+j_1+j+1}\qnf{j_1-j_{23}+j}\qnf{j_2-j_3+j_{23}}\qnf{j_2+j_3+j_{23}+1}}
{\qnf{j_{23}-j_2+j_3} \qnf{-j_1+j_{23}+j}\qnf{-j_{12}+j_3+j}\qnf{j_1+j_2+j_{12}+1}}} \times
\\ & \qquad\qquad\qquad\qquad\qquad \sqrt{\frac{\qnf{j_1+j_{23}-j}\qnf{j_2+j_3-j_{23}}
\qnf{j_1+j_2-j_{12}}\qnf{-j_1+j_2+j_{12}}}
{\qnf{-j_3+j+j_{12}}^{-1}\qnf{j_{12}+j_3+j+1}^{-1}\qnf{j_1-j_2+j_{12}}^{-1}}}\times
\\ & \sum_{l=0}^{j_1+j_2-j_{12}} \frac{(-1)^{l+j_1+j_2+j_{3}+j} [2j_{23}+2l-j_1-j_2+j_{12}+1]_q
} {\qnf{l} \qnf{j_1+j_2-j_{12}-l} \qnf{2j_{23}+l+1}\qnf{j_{23}-j_1-j_2+j_{12}+l+j_2-j_3}
\qnf{j_2+j_3-j_{23}-l}}\times \\
& \quad\dst\frac{\qnf{2j_{23}+l-j_1-j_2+j_{12}}\qnf{j_{23}+l-j_2+j_3}\qnf{-j_1+j+j_{23}+l}}
{\qnf{l+j_{12}-j_2+j+j_{23}+1}\qnf{j_1+j-j_{23}-l}\qnf{-j_1+j_{12}+j_3+j_{23}+l+1}\qnf{-j_2+j_{12}-j+j_{23}+l}},
\end{split}
\]
\end{small}
To conclude this section, let us point out that the orthogonality relations \refe{ort-rel1}
and \refe{ort-rel2} lead to the orthogonality relations for the alternative Racah polynomials
$\widetilde u_n^{\alpha,\beta}(x(s),a,b)_q$ \refe{pol-til-rac-nu-F} and their {\em duals}
$\widetilde\u_k^{\alpha',\beta'}(x(t),a',b')_q$ \refe{dua-rac-til}, respectively,
as well as the relation \refe{rel-rac-dua-rac-til} between $q$-Racah and dual $q$-Racah corresponds to
the symmetry property \refe{sim-6j}.

\subsection{Connection between $\widetilde u_{k}^{\alpha,\beta}(x(s),a,b)_q$ and
$u_{n}^{\alpha,\beta}(x(s),a,b)_q$ }

Let us obtain a formula connecting the two families $\widetilde u_{k}^{\alpha,\beta}(x(s),a,b)_q$ and
$u_{n}^{\alpha,\beta}(x(s),a,b)_q$. In fact, Eqs. \refe{6j-rac} and \refe{6j-rac-til}
sugest the following relation
between both Racah polynomials $\widetilde u_{k}^{\alpha,\beta}(x(s),a,b)_q$ and
$u_{n}^{\alpha,\beta}(x(s),a,b)_q$
\bq\label{con-u-util}
\begin{split}
& \widetilde u_{b-a-1-n}^{\alpha,\beta}(x(s),a,b)_q=
(-1)^{s-a-n}\times\\  &
 \frac{\tilG(s\!-\!a+\beta+1)\tilG(b+\alpha\!-\!s)\tilG(b+\alpha+1+s)\tilG(a+b-\beta-n)}
{\tilG(s+a\!-\!\beta+1)\tilG(\alpha+1+n)\tilG(\beta+1+n)\tilG(a+b+\alpha+1+n)}\,
 u_{n}^{\alpha,\beta}(x(s),a,b)_q.
   \end{split}
\eq

To prove it is sufficient to substitute the above formula in the difference equation
\refe{eqdif-q-new} of the $\widetilde u_{n}^{\alpha,\beta}(x(s),a,b)_q$ polynomials.
After some straightforward computations the resulting difference equation becomes into the
corresponding difference equation for the polynomials  $u_{n}^{\alpha,\beta}(x(s),a,b)_q$.

Notice that from this relation follows that
\[\begin{split}
{}_{4}\fhyp_3 & \left(\ba{c}a-b+n+1,a-b-\alpha-\beta-n, a-s,a+s+1 \\
a-b+1,2a-\beta+1, a-b-\alpha+1 \ea \,\bigg|\, q \,,\, 1 \right)\\
&
=\frac{(\beta+1|q)_{s-a} (b+\alpha+a+1|q)_{s-a}}{ (2a-\beta+1|q)_{s-a} (a-b-\alpha+1)_{s-a} }
{}_{4}\fhyp_3 \left(\ba{c} -n, \alpha+\beta+n+1, a-s,a+s+1 \\
 a-b+1,\beta+1, a+b+\alpha+1 \ea \,\bigg|\, q \,,\, 1 \right).
  \end{split}
  \]

This yield to the following identity for terminating ${}_4\phi_3$ basic series,
$n,N-n-1,k=0,1,2,\ldots$,
\[\begin{split}
{}_{4}\serie_3 & \left(\ba{c}q^{n-N+1},q^{-n-N+1}A^{-1}B^{-1}, q^{-k}, q^{-k}D \\
q^{1-N},q^{-2k}DB^{-1},q^{1-N}A^{-1} \ea \,\bigg|\, q \,,\, q \right)\\
&=\frac{q^{-kN}}{A^kB^k}
\frac{(qB;q)_k(q^{N-2k}DA;q)_k}{(q^{-2k}DB^{-1};q)_k,(q^{1-N}A^{-1};q)_k}
{}_{4}\serie_3 \left(\ba{c} q^{-n}, ABq^n, q^{-k}, q^{-k}D \\
 q^{1-N},qB,q^{N-2k}DA \ea \,\bigg|\, q \,,\, q \right).
  \end{split}
  \]
\section{Conclusions}
Here we have provided a detailed study of two kind of
Racah $q$-polynomials on the lattice $x(s)=[s]_q[s+1]_q$
and also their comparative analysis with the Racah coefficients
or $6j$-symbols of the quantum algebra $U_q(su(2))$.

To conclude this paper we will briefly discuss the relation of
these $q$-Racah polynomials with the representation theory of
the quantum algebra $U_q(su(3))$. In \cite[\S 5.5.3]{nsu} was shown that
the transformation between two different bases $(\lambda,\mu)$
of the irreducible representation of the classical (not $q$)
algebra $su(3)$ corresponding to the reductions
$su(3)\supset su(2)\times u(1)$ and $su(3)\supset u(1)\times su(2)$
of the $su(3)$ algebra in two different subalgebras $su(2)$
is given in terms of the Weyl coefficients that are, up to a
sign (phase), the Racah coefficients of the algebra $su(2)$.
The same statement can be done in the case of the quantum algebra
$su_q(3)$ \cite{as96,mal95}: The Weyl coefficients of the transformation between two
bases of the  irreducible representation  $(\lambda,\mu)$ corresponding
to the reductions $su_q(3)\supset su_q(2)\times u_q(1)$ and
$su_q(3)\supset u_q(1)\times su_q(2)$ of the quantum algebra
$su_q(3)$ in two different quantum subalgebras $su_q(2)$ coincide
(up to a sign) with the $q$-Racah coefficients of the $su_q(2)$.

In fact, the Weyl coefficients satisfy certain difference equations
that are equivalent to the differentiation formulas for the
$q$-Racah polynomials $u_n^{\alpha,\beta}(x(s),a,b)_q$ and
$\widetilde u_n^{\alpha,\beta}(x(s),a,b)_q$ so, following the
idea in \cite[\S 5.5.3]{nsu} we can assure that the main properties
of the $q$-Racah polynomials are closely related with the representations
of the quantum algebra $U_q(su(3))$. Finally, let us point out that the
same assertion can be done but with the non-compact quantum algebra
$U_q(su(2,1))$. This will be carefully done in a forthcoming paper.

\subsection*{Acknowledgements}
This research has been supported by the the DGES grant BFM 2003-06335-C03
(RAN, RCS) and PAI grant  FQM-0262 (RAN). One of authors (Yu.S.) is
thankful to the Russian Foundation of Basic research (project No 02-01-00668)
for the financial support.



{\small
\begin{thebibliography}{999}

\bibitem{renato}   R. \'Alvarez-Nodarse, { \em q-analog
of the vibrational IBM  and the
 quantum algebra $SU_q(1,1)$.} Master Thesis in Moscow State
University {\em M.V. Lomonosov}. (November 1992). (In Russian).

\bibitem{ran-tesis}  R. \'{A}lvarez-Nodarse,
{\it  Polinomios generalizados y q-polinomios: propiedades espectrales y
aplicaciones.}  {Tesis Doctoral.}
Universidad Carlos III de Madrid. Madrid, 1996. (In Spanish)

\bibitem{ran}  R. \'{A}lvarez-Nodarse,
{\it  Polinomios hipergem\'etricos y q-polinomios.} Monograf\'{\i}as del
Seminario Garc\'{\i}a Galdeano. Universidad de Zaragoza.
 Vol. {\bf 26}. Prensas Universitarias de Zaragoza, Zaragoza, Spain, 2003.
(In Spanish).

\bibitem{ran-smi}  R. \'{A}lvarez-Nodarse and Yu. F. Smirnov,
{ {\it q-}Dual Hahn polynomials on the non-uniform lattice
$x(s) = [s]_q[s+1]_q$  and the {\it q-}algebras
$SU_q(1,1)$ and $SU_q(2)$.}
{\it J. Phys. A: Math. Gen.} {\bf 29} (1996), 1435-1451.

\bibitem{as96} P. M. Asherova, Yu. F. Smirnov, and V. N. Tolstoi,
Weyl $q$-coefficients for $u_q(3)$ and Racah $q$-coefficients for $su_q(2)$,
{\it Sov. J. Nucl.Phys.} {\bf 58} (1996), 1859-1872.

\bibitem{a-w-1} R. Askey and R. Wilson,
{A set of orthogonal polynomials that generalize Racah coefficients or $6j$
symbols.}
{\it SIAM J. Math. Anal.}  {\bf 10} (1979), 1008-1020.

\bibitem{a-w-2} R. Askey and R. Wilson,
{Some basic hypergeometric orthogonal polynomials that generalize Jacobi
polynomials.}
{\it Mem. Amer. Math. Soc.}
{\bf 319}. Providence, Rhode Island, 1985.

\bibitem{arsuslov} N. M. Atakishiyev, M. Rahman, and S. K. Suslov,
{On classical orthogonal polynomials.} {\em Constr. Approx.}
{\bf 11} (1995), 181-226.

\bibitem{dri87} V.G. Drinfel'd, Quantum Groups,
{\it Proceedings of the Int. Congress of Mathematicians}.
Berkeley 1986, 798-820, American Mathematical Society. Providence,
R. I. 1987.

\bibitem{fad86} L.D. Faddev and L. A. Takhtajan,
{\it Lect. Notes Phys.} {\bf 246}, (1986) 183

\bibitem{gasper1}  G. Gasper and M. Rahman,
 {\em Basic Hypergeometric Series.} Cambridge University Press, 1990.

\bibitem{jim86} M. Jimbo, Quantum R matrix for the generalized Toda system.
{\it Commun. Math. Phys.}  {\bf 102} (1986), 537-547.

\bibitem{ki88} A. N. Kirillov and N. Yu. Reshetikhin,
{\it LOMI Preprint} E-9-88, Leningrad, 1988.

\bibitem{ks} R. Koekoek and R. F. Swarttouw,
{The Askey-scheme of hypergeometric orthogonal polynomials and its {\it q-}analogue.}
{\it Reports of the Faculty of Technical Mathematics and Informatics}
{\bf No. 98-17}. Delft University of Technology, Delft, 1998.

\bibitem{koelink} H. T. Koelink,
{Askey-Wilson polynomials and the quantum $SU(2)$ group: Survey and Applications.}
{\it Acta Applic. Math.} {\bf 44} (1996), 295-352.

\bibitem{ku81} P. P. Kulish and N. Yu. Reshetikhin,
{\it Zapiski Naucnykh Seminarov LOMI}  {\bf 101} (1981)

\bibitem{mal92}   A. A. Malashin, { \em q-analog
of Racah polynomials on the  lattice $x(s) = [s]_q[s+1]_q$ and its
connections with $6j$-symbols for the $SU_q(2)$ and $SU_q(1,1)$ quantum
algebras.}  Master Thesis. Moscow State
University {\em M.V. Lomonosov}. (January 1992). (In Russian).

\bibitem{mal95} A. A. Malashin, Yu. F. Smirnov, and  Yu.I. Kharitonov,
{\it Sov. J. Nucl.Phys.} {\bf  53 }  (1995),  1105-1119.

\bibitem{nsu-r}  A. F. \ Nikiforov, S. K. \ Suslov, and V. B. \ Uvarov,
{\it Classical Orthogonal Polynomials of a Discrete Variable.}
Nauka, Moscow, 1985, Russian Edition (in Russian).

\bibitem{nsu}  A. F. \ Nikiforov, S. K. \ Suslov, and V. B. \ Uvarov,
{\it Classical Orthogonal Polynomials of a Discrete Variable.}
{\it Springer Series in Computational Physics.}
Springer-Verlag, Berl\'{\i}n, 1991.
(Russian Edition, Nauka, Moscow, 1985)

\bibitem{nu83} A. F.  Nikiforov y V. B.  Uvarov,
{Classical orthogonal polynomials in a discrete variable
on nonuniform lattices.}
{\it Preprint Inst. Prikl. Mat. Im. M. V. Keldysha Akad. Nauk SSSR},
Mosc\'u, 1983, No. 17, (en ruso).

\bibitem{nu}  A. F. Nikiforov and  V. B. Uvarov,
{Polynomial Solutions of hypergeometric type difference Equations and their
classification}.
{\it Integral Transform.  Spec.  Funct.} {\bf 1} (1993), 223-249.

\bibitem{ro03} H. Rosengren, An elementary approach to $6j$-symbols (classical, quantum, rational and elliptic).
Preprint 2003,  arXiv:math.CA/0312310.

\bibitem{sk82} E. K. Sklyanin,
Some algebraic structures connected with the Yang-Baxter equation.
{\it Funct. Anal. Appl.} {\bf 16} (1983), 263-270.

\bibitem{smi1}  Yu. F. Smirnov, V.N. Tolstoy and Yu.I. Kharitonov, {Method of
 projection operators and the q analog of the quantum theory of angular momentum.
 Clebsch-Gordan Coefficients and irreducible
tensor operators.} {\it Sov. J. Nucl.Phys.} {\bf  53 } (1991), 593-605.

\bibitem{smi2} Yu. F. Smirnov, V.N. Tolstoy and Yu.I. Kharitonov:,
{Projection-operator
Method and the q analog of the quantum theory of angular momentum. Racah coefficients, 3j
and 6j symbols, and their symmetry properties.} {\it Sov. J. Nucl.Phys.}
{\bf  53 }  (1991),  1069-1086.

\bibitem{smi3} Yu. F. Smirnov, V.N. Tolstoy and Yu.I. Kharitonov, {Tree technique
and irreducible tensor operators for the $SU_q(2)$ quantum algebra. 9j symbols.}
{\it Sov. J. Nucl.Phys.} {\bf  55} (1992), 1599-1604 .

\bibitem{smi4} Yu. F. Smirnov, V. N. Tolstoy y Yu. I. Kharitonov,
{The tree technique and irreducible tensor operators for the $SU_q(2)$ quantum
algebra. The algebra of irreducible  tensor operators.}
{\it Physics Atom. Nucl.}  {\bf  56 } (1993), 690-700.


\bibitem{vk} N. Ja. Vilenkin and A. U. Klimyk,
{\it  Representations of Lie Groups and Special  Functions.}
{\bf  Vol. II, III}. Kluwer Academic Publishers. Dordrecht, 1992.

\end{thebibliography}

}
\end{document}